\pdfoutput=1

\documentclass[a4paper,reqno]{article}

\usepackage[colorlinks=true,hypertexnames=false,breaklinks=true,linkcolor=lightgreen,citecolor=lightgreen]{hyperref}
\usepackage{amsmath,amssymb,epsfig,float,mathabx,amsfonts,latexsym,color,algorithm,algpseudocode,xspace,fancyhdr,etoolbox}

\makeatletter
\patchcmd{\ALG@step}{\addtocounter{ALG@line}{1}}{\refstepcounter{ALG@line}}{}{}
\newcommand{\ALG@lineautorefname}{Line}
\makeatother

\newcommand{\cblue}[1]{\textcolor{lightblue}{#1}}
\definecolor{lightblue}{rgb}{0.22,0.45,0.70}
\definecolor{lightgreen}{rgb}{0.22,0.50,0.25}
\definecolor{darkred}{rgb}{0.70,0.22,0.25}
\definecolor{orange}{rgb}{1,0.5,0}
\newcommand{\cred}[1]{\textcolor{darkred}{#1}}

\usepackage{tabulary}
\newcolumntype{K}[1]{>{\centering\arraybackslash}p{#1}}

\newcommand\cero{\boldsymbol{0}}
\newcommand\OmD{\Omega^{\mathrm{D}}}
\newcommand\lambdaD{\lambda^{\mathrm{D}}}
\newcommand\lambdaE{\lambda^{\mathrm{E}}}
\newcommand\OmE{\Omega^{\mathrm{E}}}
\newcommand\muD{\mu^{\mathrm{D}}}
\newcommand\muE{\mu^{\mathrm{E}}}

\newcommand\Om{\Omega}

\newcommand{\bbeta}{\boldsymbol{\beta}}
\newcommand\bu{\boldsymbol{u}}
\newcommand\bv{\boldsymbol{v}}
\newcommand\bw{\boldsymbol{w}}
\newcommand\bz{\boldsymbol{z}}
\newcommand\bx{\boldsymbol{x}}
\newcommand\bnu{\boldsymbol{\nu}}
\newcommand\bsigma{\boldsymbol{\sigma}}
\newcommand\beps{\boldsymbol{\varepsilon}}
\newcommand\ff{\boldsymbol{f}}
\renewcommand\gg{\boldsymbol{g}}

\newcommand\bF{\mathbf{F}}

\newcommand\bG{\mathbf{G}}

\newcommand\bJ{\mathbf{J}}
\newcommand\bM{\mathbf{M}}
\newcommand\bI{\mathbf{I}}
\newcommand\PP{\mathbb{P}}
\newcommand\RR{\mathbb{R}}
\newcommand{\buE}{\boldsymbol{u}^{\mathrm{E}}}
\newcommand{\buD}{\boldsymbol{u}^{\mathrm{D}}}

\newcommand\bV{\boldsymbol{\mathcal{V}}}
\newcommand\bZ{\boldsymbol{\mathcal{Z}}}
\newcommand\cD{\mathcal{D}}
\newcommand\cQ{\mathcal{Q}}
\newcommand\cR{\mathcal{R}}
\newcommand\cS{\mathcal{S}}

\newcommand{\cT}{\mathcal{T}}

\DeclareMathAlphabet{\mathbfsf}{\encodingdefault}{\sfdefault}{bx}{n}
\newcommand{\bsA}{\mathbfsf{A}}
\newcommand{\bsB}{\mathbfsf{B}}
\newcommand{\bsC}{\mathbfsf{C}}
\newcommand{\bsCu}{\mathbfsf{C}_u}
\newcommand{\bsD}{\mathbfsf{D}}
\newcommand{\bsJ}{\mathbfsf{J}}
\newcommand{\bsM}{\mathbfsf{M}}

\newcommand{\bdiv}{\operatorname*{\mathbf{div}}}
\newcommand{\vdiv}{\operatorname*{div}}
\newcommand{\pt}{\partial_t}
\newcommand{\dx}{\,\mbox{d}\bx}
\newcommand{\ds}{\,\mbox{d}s}

\newcommand{\dtnum}[1]{\Delta t^{#1}}
\newcommand{\tn}[1]{t_{#1}}
\renewcommand{\ng}{n+\gamma}

\newcommand{\vstep}[3]{\mathbf{#1}_{#2}^{#3}}
\newcommand{\dvstep}[5]{\delta_{#1}^{#2}\vstep{#3}{#4}{#5}}
\newcommand{\ewo}{\boldsymbol{e}_w^1}
\newcommand{\ewt}{\boldsymbol{e}_w^2}

\newcommand{\MJ}{\frac{1}{\dtnum{}}\bsM-\gamma\bsJ}
\newcommand{\MJs}{\frac{1}{\dtnum{}}\bsM^*-\gamma\bsJ^*}
\newcommand{\res}[1]{\mathrm{res}_{\mathrm{#1}}}
\newcommand{\eN}{e(\dvstep{k+1}{}{w}{z}{},\widehat{sc})}

\newcommand{\True}{\textbf{true}\xspace}
\newcommand{\False}{\textbf{false}\xspace}

\renewcommand{\And}{\textbf{and}\xspace}

\newcommand{\Break}{\textbf{break}}

\numberwithin{equation}{section}
\numberwithin{table}{section}
\numberwithin{figure}{section}
\allowdisplaybreaks

\newcommand{\executeiffilenewer}[3]{%
\ifnum\pdfstrcmp
{\pdffilemoddate{#1}}%
{\pdffilemoddate{#2}}%
>0%
{\immediate\write18{#3}}%
\fi%
}
\newcommand{\includeeps}[1]{%
\IfFileExists{#1.pdf}%
{\executeiffilenewer{#1.eps}{#1.pdf}{epspdf #1.eps}}%
{\immediate\write18{epspdf #1.eps}}%
}

\title{Numerical approximation of a mechanochemical \\interface  
model for skin appendage\thanks{This work was partially supported by grants 
from the Swiss National Science Foundation (FNSNF, grants 31003A\_140785 
and SINERGIA CRSII3\_132430), by the SystemsX.ch initiative (project Epi-PhysX); 
and by the Engineering and Physical Sciences Research Council (EPSRC, research grant EP/R00207X/1).}}

\author{\normalsize {\sc Luis Miguel De Oliveira Vilaca}\thanks{Laboratory of Artificial \& Natural Evolution (LANE), 
Department of Genetics and Evolution, University of Geneva, 4 Boulevard d'Yvoy,  
1205 Geneva, Switzerland; and SIB Swiss Institute of Bioinformatics, Geneva, Switzerland. 
Email: {\tt \{LuisMiguel.DeOliveiraVilaca,Michel.Milinkovitch@unige.ch\}@unige.ch}.},  \ 
{\sc Michel C. Milinkovitch$^{\dag}$},  \ 
{\sc Ricardo Ruiz-Baier}\thanks{Corresponding author. Mathematical Institute, Oxford University, Andrew Wiles Building, 
Woodstock Road, OX2 6GG Oxford, UK. Email: {\tt RuizBaier@maths.ox.ac.uk}.}}

\date{\today}
 
\begin{document}
\maketitle

\begin{abstract}
\noindent We introduce a model for the mass transfer of molecular activators and inhibitors in two media separated by an interface, and study its interaction with the deformations exhibited by the two-layer skin tissue where they occur. The mathematical model  results in a system of nonlinear advection-diffusion-reaction equations including cross-diffusion, and coupled with  an interface elasticity problem. We propose a Galerkin method for the discretisation of the set of governing equations, involving also a suitable Newton linearisation, partitioned techniques, non-overlapping Schwarz alternating schemes, and high-order adaptive time stepping algorithms. The experimental accuracy and robustness of the proposed  partitioned numerical methods is assessed, and some illustrating tests in 2D and 3D are  provided to exemplify the coupling effects between the mechanical properties and the advection-diffusion-reaction interactions involving the two separate layers. 
\end{abstract}



\section{Introduction and modelling considerations}


Reliable modelling of mechanochemical properties of composite materials is of key importance in a variety of engineering, material science and life science applications. These include tissue engineering, wound healing manipulation, bone fracture repair, design of photovoltaic devices, polymer adhesion, and many others. Of special interest to us are the two-way interactions between, on one hand, biological tissue deformation and, on the other hand, migration and proliferation of certain cell types during organismal development.

One application of such a general model is the study of the vertebrate skin appendage development and evolution, especially their diversity of forms (such as hairs, spines, feathers, and scales;  \cite{glover2017,MILLAR2002,montandon14,Painter433}), and spatial organisations. Note that we have recently shown that all skin appendages are homologous structures developing from an anatomical placode and associated signalling molecules \cite{dipoi2016}.

Turing proposed in his seminal work \cite{turing} that pattern formation in biological systems could be explained by a system of diffusion and interacting chemicals that he called 'morphogens'. Many subsequent studies \cite{liu2006,Painter433,shoji,Sick1447} investigated biological patterning in the framework of such purely chemical reaction-diffusion models. However, as already suggested by Turing \cite{turing}, mechanical aspects might also be relevant to biological patterning processes. During the last decade, biological experiments have clearly indicated that not only chemical (pattern of morphogen concentrations), but also mechanical (stresses and strains), parameters are key to skin development \cite{ahtiainen14,AMAN2010,glover2017,Oliver94,TRANQUI2000} and are synchronised during embryonic development \cite{MURRAY1988}.

Although the full detailed underlying mechanisms are still incompletely understood, authors have attempted to integrate both chemical and mechanical parameters in their mathematical models of skin patterning \cite{MurrayBookMB2SMBA,OsterMAMM}. The various proposed mechanochemical models differ by the nature of the interactions between morphogen concentration and mechanics. For instance, while in \cite{MurrayGBPF,quarteroni17} mechanical responses are directly triggered by chemical signalling, other models \cite{moreo10,MurrayGBPF} propose that migration of cells directly generate local forces on the tissue. Such interactions between cells and morphogens can be complemented with chemical and mechanical exchanges between the two layers (dermis and epidermis) that constitutes the skin \cite{cruywagen1992,shawAMCSP}. In one example \cite{badugu12}, a three-chemical species model of limb bud development integrates growth velocity of the tissue normal to the limb surface and that directly depends on the local concentration of fibroblast growth factor. Other models consider the tissue growth as a liquid displacement \cite{chertock12}, or aggregation of cells is the result of a chemical pre-pattern.

Here we focus on the computational modelling of mechanochemical interactions (inspired by \cite{cruywagen1992} and \cite{shawAMCSP}) based on simple reaction-diffusion systems governing the relations between, on one hand, the concentration of morphogenic proteins and cells densities and, on the other hand, mechanical stresses and strains. The mathematical description proposed in these two references is sufficiently general to be applied to multiple biological contexts (e.g., angiogenesis, fingerprint formation, cartilage condensation \cite{MurrayBookMB2SMBA}), justifying the numerical algorithm presented here. Experimental evidence shows that the skin short-term and low-amplitude response to mechanical stimuli is that of an elastic solid \cite{dillon99}. We therefore assume that the two-layer domain of the skin (dermis and epidermis) can deform according to its inherent mechanical response, and that it is affected isotropically via an internal force that depends on the signalling factors. In addition, given the long spatio-temporal scale of the organismal developmental process, short-term fluctuations and inertial effects can probably be safely discarded to simplify the mathematical model \cite{moreo10}. We therefore adopt the two following fundamental assumptions: (i) the equilibrium of forces in the system is established by a quasi-static balance of linear momentum, and (ii) chemical species, governed by an advection-diffusion reaction (ADR) system, affect the medium deformation (through coupling functions). The chemical species' dynamics represents the spatio-temporal evolution of two morphogenic proteins and/or of cells that trigger and control the internal forces affecting the momentum balance. In turn, cell motion induces local and global traction forces on the elastic body, leading to substrate deformations \cite{Oliver94}. We point out that cell motility depends on different factors that include kineses (increase in motile activity without any directional component), taxes (directed motion up or down a gradient), and guidances (motion directed by substratum cues, also called contact guidance) \cite{MurrayGBPF} that can also act locally or globally. Although long-range effects can be implemented through high-order spatial derivatives (appearing in the ADR or force conservation equations), we will here incorporate only short-range effects because our focus is on the numerical and algorithmic aspects of interface problems (such as the dermis-epidermis interactions in the skin).

As suggested elsewhere \cite{cruywagen1992}, here we distinguish between the two very different cellular dynamics of the dermis and the epidermis. Whereas we adopted a general logistic growth mitotic rate and non-negligible short range cell diffusion \cite{MURRAY1988} in the dermis, we account for a zero production rate of cells and very small diffusivities in the equations governing epidermal variables \cite{ahtiainen14,shawAMCSP}. This difference of treatment in the two tissues is justified by their respective structures: cells in epidermal tissue form sheets of adjoined interacting cells, while mesenchymal cells in the dermis are more mobile and separated by substantial extracellular matrix. However, as epidermal cells are not completely attached to each other and can rearrange topologically \cite{ahtiainen14}, we assume in both dermis and epidermis that all diagonal terms in the diffusivity matrices are strictly positive.

Regarding the interface problem between dermis and epidermis, we considered multiple methods for numerical approximation of linear elasticity using mixed finite elements: parallel and sequential Schwarz iterations using various interface conditions based on boundary element-finite couplings \cite{gatica97}, $hp$-mortar \cite{belgacem03}, Nitsche \cite{becker09}, multiscale settings \cite{buck13}, or matching interface algorithms \cite{wang15}.
We finally chose a Schwarz iteration method with domain decomposition where the coupled set of mechanochemical equations is solved separately on each side of the interface, imposing mixed-type (Robin) boundary conditions, i.e., both fluxes and density values are equal in both tissues at their interface. A rich body of literature deals with choosing the weighting coefficients on the Robin interface conditions in order to optimise the convergence of the alternating Schwarz algorithms under different scenarios (see e.g. \cite{chen14,gander02,gosselet07,qin08} and references therein). As the optimal weighting scheme tends to be problem-dependent and might change substantially for diffusion, elasticity, or fluid-oriented problems, we will only describe general schemes and mention which specific weighting best suits the problems we investigate here. 

Here, we define a family of finite element methods for the semi-discretisation of the problem, where the classical MINI-element \cite{arnold84} is used for the approximation of both displacement and pressure on each subdomain, whereas piecewise linear and continuous elements are used to discretise the vector of species concentration. Concerning the fully-discrete (spatial and time discretisation) scheme for the model problem, we apply an adaptive time-stepping strategy based on diagonally-implicit Runge-Kutta methods. This so-called TR-BDF2 method was introduced in the context of electronic circuit simulation \cite{bankTRBDF2} (see also \cite{bonaventura17,HOSEA1996,edwards11} for further details) and has the advantage of using the same Newton iteration matrix to solve different temporal stages within a time step, a particularly interesting feature when addressing large systems of ODEs. The step size adaptivity exploits classical automatic control techniques used in the integration of so-called stiff systems \cite{HairerBookODEI,HairerBookODEII}, i.e., problems that tend to be numerically unstable.

This article is organised in the following manner. The remainder of this section contains preliminary notation and generalities to be used throughout the paper. Section~\ref{sec:model} provides a detailed description of the differential equations constituting the model problem, with special emphasis on the coupling occurring at the interface. In Section~\ref{sec:FE}, we introduce the Galerkin method, we present the adaptive time discretisation, and we specify algorithmic considerations. Section~\ref{sec:results} contains a collection of computational tests focusing on the accuracy of the numerical schemes and on the suitability of the model equations to represent mechanochemical problems. We close with a summary and discuss potential extensions of the model in Section~\ref{sec:concl}.

Throughout the paper, we let $\Om\subset\RR^d$, $d\in\{2,3\}$ 
denote a deformable body with polyhedral boundary $\partial\Om$, and denote by
$\bnu$ the outward unit normal vector on $\partial\Om$.  
Moreover, given a generic domain $\cD\in \RR^d$, $d\in\{2,3\}$, for 
two vectorial functions $\bu,\bv\in L^2(\cD)^d$ 
will denote their inner product as $(\bu,\bv)_{\cD} =\int_{\cD} \bu\cdot\bv\dx$, whereas 
if $\cS$ is a $(d-1)-$dimensional surface then their pairing will be written as   
$\langle \bu,\bv\rangle_{\cS} = \int_{\cS} \bu\cdot\bv\ds$.   In addition, by 
$\PP_r(\cD)$ we will denote the space of polynomial functions of degree $s\leq r$ 
defined on the domain $\cD$, and $\bI$ stands for the identity matrix in $\mathbb{R}^{d\times d}$. 
Finally,  we will employ $\cero$ to denote a
generic null vector or the zero operator.

\section{Set of governing equations}\label{sec:model}
\begin{figure}[!t]
\begin{center}
\begin{picture}(100,150)
\put(0,0){\includegraphics[width=0.4\textwidth]{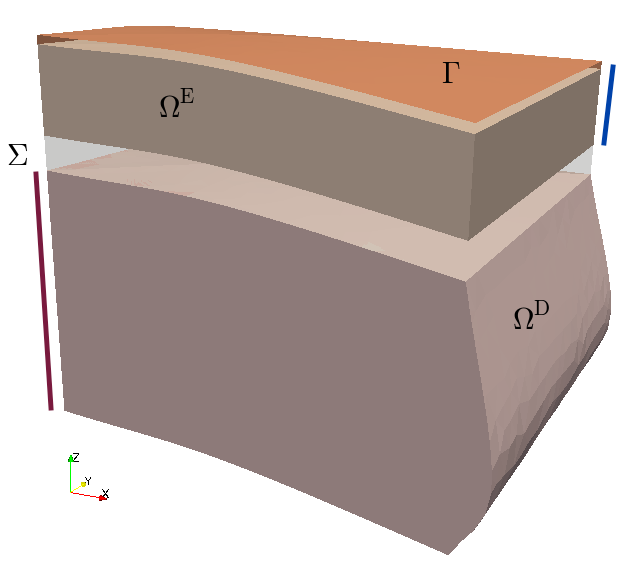}}
\put(-10,60){\cred{$\partial\OmD$}}
\put(175,125){\cblue{$\partial\OmE$}}
\end{picture}
\caption{Schematic representation of dermis and epidermis layers separated by the interface $\Sigma$ and 
an exposed epidermic surface is denoted by $\Gamma$.}\label{fig:ex00sketch} 
\end{center}
\end{figure}

\paragraph{Linear elastic solids.} 
We assume that  the domain $\Om$ is disjointly 
divided into two connected subdomains (with different material properties) denoted 
by $\OmD$ and $\OmE$, respectively, whose intersection constitutes the 
interface between the regions and is denoted $\Sigma=\partial\OmD\cap \partial\OmE$. 
On $\partial\OmE$ we define a part of the boundary related to the exposed surface, denoted as $\Gamma$ 
(see Figure~\ref{fig:ex00sketch}). 

This heterogeneity of the skin (corresponding to the presence of an homogeneous \emph{dermis} and an homogeneous \emph{epidermis}, occupied by $\OmD$ and $\OmE$, respectively) is accounted for by the spatially-dependent (but locally constant) Lam\'e coefficients  
$$ 
\lambda = \begin{cases} 
\lambdaD & \text{in $\OmD$},\\
\lambdaE & \text{in $\OmE$},\end{cases}
\quad \mu = \begin{cases} 
\muD & \text{in $\OmD$},\\
\muE & \text{in $\OmE$}.\end{cases}
$$
We also suppose that the two-layered elastic medium is nearly incompressible, 
so, in view of avoiding the well-known phenomenon of volumetric locking 
(inability of some discrete spaces to approximate the true displacement while 
satisfying incompressibility, see e.g. \cite{belgacem03}), 
a mixed formulation is adopted in terms of displacements $\bu^*$ 
and the pressure $p^* = -\lambda^* \vdiv\bu^*$, with $*\in\{$D,E$\}$ 
(although the interpretation of $p$ as pressure holds exclusively in 
the incompressibility limit where $\lambda^*/\mu^*$ takes on quite large 
values). For all $0<t<T$, and given a body force $\ff$ (whose specific form depends on the vector $\bw$; see equation \eqref{2way-coupling} below) applied on each domain, we resolve the displacement / pressure pair $(\bu^*(t),p^*(t)):\Om^*\to [\RR^d \times \RR]$ on each subdomain, satisfying spring conditions at the surface boundary (Robin conditions), having clamped boundaries elsewhere, and assuming an  adhesion 
condition at the interface (matching displacements and preservation of traction forces imply continuity of the medium across $\Sigma$). Then the global linear Navier-Lam\'e problem adopts the form 
\begin{align}
\nonumber
-\bdiv \bigl(2\mu^* \beps(\bu^*) - p^* \bI \bigr) & = \ff(\bw^*)  & \text{in $\Om^*\times(0,T)$}, \\
p^* + \lambda^* \vdiv\bu^*& =0  &\text{in $\Om^*\times(0,T)$}, \nonumber\\
\bsigma^\mathrm{E}\bnu + \alpha^\mathrm{E} \buE & = \cero  &\text{on $\Gamma\times(0,T)$},\label{eq:elasticity} \\
 \bu^* & = \cero  &\text{on $\partial\Om\setminus\Gamma\times(0,T)$}, \nonumber\\
\bsigma^{\mathrm{D}}\bnu = \bsigma^{\mathrm{E}}\bnu,& \quad  \buD   = \buE  &\text{on $\Sigma\times(0,T)$},\nonumber 
\end{align}
where $\beps(\bu^*) = \frac{1}{2}(\nabla\bu^*+\nabla\bu^{*T})$ is the tensor of 
infinitesimal strains (symmetric gradient of displacement), $\bsigma^*=2\mu^* \beps(\bu^*) - p^* \bI $ is the Cauchy stress tensor, and $\alpha^*>0$ is the Robin coefficient. The extension to non-homogeneous displacement conditions ($\bu|_{\partial\Om} \neq 0$) can be handled in a standard way by adding suitable liftings on the data. Below, we use the asterisk superscript $*\in\{$D,E$\}$ to denote quantities related to each subdomain, and we adopt the convention that the unit normal vector $\bnu$ on the interface is pointing from $\OmD$ to $\OmE$. 

We stress that modelling the development of soft tissue and cells, in terms 
of mechanical dynamic behaviour is still an active area of research \cite{allena2013,Wyscz2012}. 
Biological tissues, depending on the nature of the constituent cells, present a large variability in their elastic properties (including viscous, nonlinear, or nonlocal effects) which can change substantially through 
spatio-temporal scales \cite{merker2016,Pawlaczyk2013,Wyscz2012}. Some recent analyses were based on nonlinear St. Venant-Kirchhoff descriptions \cite{merker2016}, including viscoelastic terms \cite{allena2013,MurrayGBPF,Tosin2006,turing}, and even fluid dynamics, as mentioned in \cite{Wyscz2012}. A general consensus on how to represent mechanical properties of the skin tissue under morphogenesis is not yet available. Hence, our description will be restricted to classical linear elastic materials because it accurately represents small tissue deformations manifested during the first stage of skin appendage formation. 

\paragraph{Advection-diffusion-reaction equations on an underlying deforming medium.}  
Let us focus on the spatio-temporal interaction
of the densities of $m$ species $\bw = (w_1, \ldots ,w_m)$. We assume that 
all constituent species are equi-present at each spatial point (i.e., they are spatially continuous), so applying the Reynolds transport theorem to their mass conservation (assuming zero-flux boundary conditions at the surface $\Gamma$ and continuity of fluxes and densities at the interface) leads to the system 
\begin{align}
\pt \bw^* + (\pt\bu^*\cdot\nabla) \bw^* - \bdiv (\bM^* \nabla\bw^*) 
 & = \bG(\bw^*) + \gg(\bu^*) & \text{in $\Om^*\times(0,T)$},  \nonumber\\
 \bM^* \nabla\bw^* \bnu & = \cero & \text{on $\partial\Om \times(0,T)$},  \label{eq:rd}\\
\bM^\mathrm{D} \nabla\bw^\mathrm{D} \bnu  = \bM^\mathrm{E} \nabla\bw^\mathrm{E} \bnu, & \quad  
\bw^\mathrm{D} = \bw^\mathrm{E} & \text{on $\Sigma \times(0,T)$}.\nonumber 
\end{align} 
Here $\bM^*\in\RR^{m\times m}$ is a (not necessarily diagonal, but positive definite) 
tensor of self- and cross-diffusion rates, $\bG:\Om\to\RR^m$ is a vectorial function 
containing the reaction kinetics of the system and representing the 
production and degradation of species concentration on each subdomain $\OmD,\OmE$; 
and the term $\gg$ encodes an additional coupling contribution (detailed in equation \eqref{2way-coupling} below). The system is endowed with adequate initial conditions for $\bw^*$, and the general form \eqref{eq:rd} 
accommodates either dimensional or dimensionless descriptions. 

Here,we focus on two levels of model complexity. Firstly, we define a simple 
two-variable activator-inhibitor model with quadratic reaction terms due to Gierer-Meinhardt 
(cf. \cite{gierer72}): 
\begin{equation}\label{reaction:GM}
\bG(\bw^*) = \begin{pmatrix} 
\rho_2(\rho_0+\rho_1 \frac{(w_1^*)^2}{w^*_2}) - \rho_3 w^*_1 \\
\rho_4 (w^*_1)^2 - \rho_5w^*_2\end{pmatrix},\end{equation} 
where  $w^*_1$ is the short-range autocatalytic substance 
 and $w^*_2$ is the long-range antagonist; and the $\rho_i$'s are positive model parameters.
 
Second, we introduce a more general model of skin appendages development were $w^*_1$ is  the cell density, 
$w^*_2$ is the density of the extracellular matrix, and $w^*_3,w_4^*$ are the concentrations of two 
morphogen signalling chemical compounds. All quantities are assumed to diffuse, so none 
of the diagonal terms in $\bM^*$ is zero (nevertheless, self-diffusion for cell concentration in the epidermis and the matrix density in the two tissues, will be very small in comparison to the other values in the diagonal of $\bM^\mathrm{E}$), and we discard higher order terms (derivative with order higher than two) representing haptotaxis and non-local effects in the cell migration fluxes. The reactions 
consist in a logistic relation for the mitotic rate, a zero rate of matrix material secretion, 
and the production, as well as degradation, of morphogens:  
\begin{equation}\label{reaction:general}
\bG(\bw^\mathrm{D}) = \begin{pmatrix}
r_1 w_1^\mathrm{D}(r_0 - w_1^\mathrm{D})\\
0 \\
r_2 w_1^\mathrm{D} - r_3 w_3^\mathrm{D} \\
- r_4  w_1^\mathrm{D} w_4^\mathrm{D} 
\end{pmatrix}, \quad 
\bG(\bw^\mathrm{E}) = \begin{pmatrix}
0\\
0\\
- r_5 w_1^\mathrm{E}  w_3^\mathrm{E} \\
 r_6 w_1^\mathrm{E} - r_7 w_4^\mathrm{E} 
\end{pmatrix},
\end{equation}
where $r_0$ is the maximum cell density, and $r_i$'s are positive model constants. The specifications in 
\eqref{reaction:general} have been adapted from the skin pattern formation model proposed in \cite{cruywagen1992}. Similar models can be found in the context of cell migration, angiogenesis, follicle differentiation, 
or mesenchymal morphogenesis (see e.g.~\cite{ahtiainen14,manoussaki03,montandon14,OsterMAMM,vaughan13}). 

\paragraph{Coupling mechanisms.} Extensive experimental evidence 
indicates that both biochemical and mechanical processes work in concert during pattern formation 
(see for instance \cite{ahtiainen14,allena2013,glover2017,lecuit2007,merker2016,Shyer811}). 
In the specific context of skin appendage formation, signalling pathways control cell movements that eventually lead to mechanical tissue deformation due to cell aggregation \cite{ahtiainen14,glover2017,Shyer811}. These interactions can be represented by diverse 
local or non-local processes including advection, cytoskeleton regulation, chemical gradients, 
active stress, strain dependent diffusivity, and many others \cite{MURRAY1988,moreo10}. In the present 
model we opt for including only local variations of the ADR and mechanical variables. 
This simplification rules out some phenomena such as inter-cellular force transmission through cytoskeleton,  
nonlinear interactions (e.g., cells sensing, through filopodia, long-range variations), 
or anisotropy effects, which can be of special importance in some biological contexts 
\cite{cruywagen1992,manoussaki03}. Nonetheless, as suggested in \cite{moreo10}, local models such as the one we use here, may suffice to reproduce accurately a variety of mechanical processes without the need for additional constitutive relations. 

The PDE-based model \eqref{eq:elasticity}-\eqref{eq:rd} (written in Eulerian form and solving for 
displacements, solid pressure, and molecular variables)  assumes that changes in the chemical 
concentrations do not affect the mechanical properties of the solid, nor the Lam\'e constants; 
the species only act as forcing terms on the momentum balance, i.e., the intensity of the deformation will depend locally on the gradient of the species concentrations. Conversely, we here suppose
that the body deformation affects the species dynamics by means of advection and through the term $\gg(\bu^*)$ on the right-hand side in \eqref{eq:rd}, carrying local information about the medium dilation. Therefore, we set 
\begin{equation}\label{2way-coupling}
 \ff(\bw^*) = c_f^* \sum_{i=1}^m \nabla w^*_i, \qquad \gg(\bu^*) = c_g^*\vdiv\bu^*\begin{pmatrix}1\\ \vdots \\1\end{pmatrix},
 \end{equation}
where $c_f^*,c_g^*$ are  model parameters. Note that the form 
of $\ff$ is equivalent to considering a stress component depending on the 
concentrations of some $w_i^*$, that is $\bsigma^* = \bsigma^*_{\mathrm{eff}} + w^*_i\bI$, 
where $\bsigma^*_{\mathrm{eff}}$ is an effective elastic stress. Similar coupling conditions can be found  
in \cite{gatica18,lenarda17,ruiz15}.

\subsection{Weak formulation for the multidomain problem}
Multiplying \eqref{eq:elasticity}-\eqref{eq:rd} by suitable test functions, 
integrating by parts over each region $\Om^*$, and exploiting the boundary 
and interface conditions,  we arrive  
at the following primal-mixed variational problem associated to the elasticity-ADR 
problem: For $t>0$, find  $(\bu^*(t),p^*(t),\bw^*(t)) \in \bV^*\times \cQ^*\times\bZ^*$ such that 
\begin{align}
a(\bu^*(t),\bv)   + b(p^*(t),\bv) + \alpha^*\langle \bu^*(t),\bv\rangle_\Gamma & = F(\bw^*; \bv) & \forall \bv \in \bV^*,\nonumber\\
b\bigl(q,\bu^*(t)\bigl) - c(p^*(t),q) & = 0 &\forall q\in \cQ^*,\label{weak-nonlinear}\\
\bigl(\pt\bw^*(t), \bz\bigr)_{\Om^*} + C(\bu^*(t),\bw^*(t);\bz) + d(\bw^*(t),\bz) & = G(\bw^*;\bz) &\forall \bz\in\bZ^*, \nonumber
\end{align}
where $\alpha^\mathrm{E}$ is the coefficient associated to the Robin boundary condition 
at the surface, and the specific boundary treatment implies that adequate functional spaces are 
$$\bV^*   : = H^1_0(\Om^*)^d = \{ \bv\in H^1(\Om^*)^d: \bv|_{\partial\Om\setminus\Gamma} = \cero\},\qquad 
\cQ^* = L^2(\Om^*), \qquad \bZ^* = H^1(\Om^*)^m.$$
Linear and nonlinear forms (lowercase and uppercase, respectively) involved in \eqref{weak-nonlinear} are defined as 
\begin{align*}
a(\bu,\bv) &: = 2\mu^* ( \beps(\bu),\beps(\bv))_{\Om^*}, \quad  
b(q,\bv) : = - (q ,\vdiv \bv)_{\Om^*}, \quad c(p,q):=\frac{1}{\lambda^*}(p,q)_{\Om^*},\\ 
C(\bu(t),\bw;\bz) & := \bigl([\pt\bu\cdot\nabla] \bw,\bz\bigr)_{\Om^*}, \quad 
d(\bw,\bz) : =  (\bM^* \nabla\bw,\nabla\bz)_{\Om^*},\\
F(\bw;\bv) & : = \bigl(\ff(\bw),\bv \bigr)_{\Om^*}, \quad G(\bu,\bw;\bz): = \bigl(\bG(\bw)+g(\bu),\bz\bigr)_{\Om^*}.
\end{align*}
The unique solvability of \eqref{weak-nonlinear} has not yet been addressed in the literature, 
however related 
results concerning the coupling of linearised mechanics and reaction-diffusion systems have been 
recently addressed in the context of cardiac and plant cell biomechanics \cite{andreianov15,ptashnyk16}.

\paragraph{Treating the interface: a non-overlapping continuous Schwarz method.} As mentioned previously, the treatment of the interface (here, between the epidermis and dermis) is an important issue for interacting media problems and this paragraph will focus on its realisation. Regarding the geometric separation between the problems defined on each side of the interface, we remark that \eqref{weak-nonlinear} is valid in particular if we consider a solution globally defined in $\Om$. Nevertheless, in the discrete case we aim at solving 
two separate problems on each subdomain and therefore we require to impose jump conditions  in an adequate manner.  We adopt a  non-overlapping Schwarz iteration based on Robin-Robin type interface conditions for the ADR system
\begin{equation}\label{eq:interface-adr}
\cR(\bw^\mathrm{D}):=\bM^\mathrm{D} \nabla\bw^\mathrm{D} \bnu + K^\mathrm{D}\bw^\mathrm{D} = \bM^\mathrm{E} \nabla\bw^\mathrm{E} \bnu + K^\mathrm{E}\bw^\mathrm{E}=:\cR(\bw^\mathrm{E})  \qquad \text{on $\Sigma \times(0,T)$},
\end{equation}
where $\cR(\cdot)$ denotes the Robin transmission operator, and 
$K^\mathrm{D},K^\mathrm{E}$ are nonnegative acceleration 
constants satisfying $K^\mathrm{D}+K^\mathrm{E}>0$  (see, for instance \cite{lions90} and \cite{dolean15,quarteroni99}). On the other hand, the adhesion 
interface condition featured by the elasticity problem is incorporated through the following relation 
\begin{equation}\label{eq:interface-elast}
\cS(\buD,p^\mathrm{D}):= \bsigma^{\mathrm{D}}\bnu +J^\mathrm{D}\buD = \bsigma^{\mathrm{E}}\bnu + J^\mathrm{E}\buE =: \cS(\buE,p^\mathrm{E}) \qquad 
\text{on $\Sigma\times(0,T)$},
\end{equation}
where $\cS(\cdot,\cdot)$ is a transmission operator associated to the elasticity problem, and $J^*> 0 $ is 
another acceleration constant. The classical iterative Schwarz algorithm 
(still continuous and still nonlinear) then reads: 
For $t>0$ and for $\ell = 0,\ldots $, find  $(\bu^{*,\ell+1},p^{*,\ell+1},\bw^{*,\ell+1}) \in \bV^*\times \cQ^*\times\bZ^*$ such that 
\begin{align}
a(\bu^{*,\ell+1},\bv)   + b(p^{*,\ell+1},\bv) + L^*\langle \bu^{*,\ell+1},\bv\rangle_{\Sigma\cup \Gamma} & = F(\bw^{*,\ell+1}; \bv)+ \langle\cS^{\dag,\ell},\bv\rangle_\Sigma & \forall \bv\in \bV^*,\nonumber\\
b\bigl(q,\bu^{*,\ell+1}\bigl) - c(p^{*,\ell+1},q) & = 0 &\forall q\in \cQ^*,\label{weak-schwarz}\\
\bigl(\pt\bw^{*,\ell+1}, \bz\bigr)_{\Om^*} + C(\bu^{*,\ell+1},\bw^{*,\ell+1};\bz) \qquad \nonumber \\ 
+d(\bw^{*,\ell+1},\bz) + K^*\langle \bw^{*,\ell+1},\bz\rangle_\Sigma & = G(\bw^{*,\ell+1};\bz) +\langle\cR(\bw^{\dag,\ell}),\bz\rangle_\Sigma  &\forall \bz\in\bZ^*, \nonumber
\end{align}
where 
$$L^* =\begin{cases}
\alpha^\mathrm{E} & \text{if $* = $\,E and on $\Gamma$},\\
J^*& \text{ on $\Sigma$},\\
0 & \text{otherwise},\end{cases}$$
and the superscript $\dag\in\{$D,E$\}$ denotes quantities in the subdomain 
opposite to $*$. Here we have dropped the explicit time-dependence whenever clear from the context.

\section{A finite element method and its splitting scheme}\label{sec:FE}
In this section we detail a Rothe-type strategy employed to numerically solve the set of governing 
equations \eqref{eq:elasticity}-\eqref{eq:rd}. We first apply a discretisation in space written as a 
primal-mixed Galerkin formulation, which we later combine with fixed-point schemes and a Schwarz iterative algorithm 
decoupling the problem into four main blocks corresponding to each side of the interface, and to the elasticity and 
 ADR equations. Then we discuss the 
time discretisation based on an adaptive scheme and operator splitting. 

\subsection{Nonlinear, partitioned Galerkin formulation} 
Let us denote by $\cT_h^*$ a shape-regular simplicial partition of $\Om^*$ by 
triangles (for $d=2$, or tetrahedra for $d=3$) $K$ of diameter 
$h_K$ and define the global meshsize as $h:=\max\{h_K: K\in \cT^*_h\}$. For 
clarity of the presentation, we assume that the two meshes (associated with each subdomain) match at the interface.
We will seek the approximation $(\bu_h^*,p_h^*,\bw_h^*)$ of 
$(\bu^*,p^*,\bw^*)$, in the finite dimensional product space 
\begin{equation}\label{eq:subspaces}
[\bV^*_h\times \cQ^*_h \times \bZ^*_h] \subset [\bV^*\times \cQ^*\times \bZ^*],
\end{equation}
to be made precise below. Then, 
the semidiscrete Galerkin scheme associated to the Schwarz 
alternating algorithm \eqref{weak-schwarz} and a fixed-point decoupling of elasticity and 
ADR subproblems can be recast in the following partitioned form, iterated until 
convergence: 
For a given  $t>0$,  and given $(\bu_h^{*,\ell},p_h^{*,\ell},\bw_h^{*,\ell}) \in \bV_h^*\times \cQ_h^*\times\bZ_h^*$, 
find the next iteration $(\bu_h^{*,\ell+1},p_h^{*,\ell+1},\bw_h^{*,\ell+1})$ from the sequential 
system 
\begin{align}
&\begin{cases}
& \bigl(\pt\bw_h^{\mathrm{D},\ell+1}, \bz_h^\mathrm{D}\bigr)_{\Om^\mathrm{D}} + C(\bu_h^{\mathrm{D},\ell},\bw_h^{\mathrm{D},\ell+1};\bz_h^\mathrm{D}) 
+d(\bw_h^{\mathrm{D},\ell+1},\bz_h) \\
& \qquad + K^\mathrm{D}\langle \bw_h^{\mathrm{D},\ell+1},\bz_h^\mathrm{D}\rangle_\Sigma - G(\bw_h^{\mathrm{D},\ell+1};\bz_h^\mathrm{D})  = \langle\cR_h^{\mathrm{E},\ell},\bz_h^\mathrm{D}\rangle_\Sigma   \quad \forall \bz_h^\mathrm{D}\in\bZ_h^\mathrm{D},
\end{cases}\label{galerkin-wD}\\[2ex]
&\begin{cases}
 a(\bu_h^{\mathrm{D},\ell+1},\bv_h^\mathrm{D})   + b(p_h^{\mathrm{D},\ell+1},\bv_h^\mathrm{D}) + L^\mathrm{D}\langle \bu_h^{\mathrm{D},\ell+1},\bv_h^\mathrm{D}\rangle_{\Sigma\cup \Gamma} \!\!\! & = F(\bw_h^{\mathrm{D},\ell}; \bv_h^\mathrm{D})+ \langle\cS_h^{\mathrm{E},\ell},\bv_h^\mathrm{D}\rangle_\Sigma \ \forall \bv_h^\mathrm{D}\in \bV_h^\mathrm{D},\\
 b\bigl(q_h^\mathrm{D},\bu_h^{\mathrm{D},\ell+1}\bigl) -\, c(p_h^{\mathrm{D},\ell+1},q_h^\mathrm{D})  & = 0 \quad \forall q_h^\mathrm{D}\in \cQ_h^\mathrm{D},
\end{cases}\label{galerkin-upD}\\[2ex]
&\begin{cases}
& \bigl(\pt\bw_h^{\mathrm{E},\ell+1}, \bz_h\mathrm{E}\bigr)_{\Om^\mathrm{E}} + C(\bu_h^{\mathrm{E},\ell},\bw_h^{\mathrm{E},\ell+1};\bz_h^\mathrm{E})  
 +d(\bw_h^{\mathrm{E},\ell+1},\bz_h^\mathrm{E}) \\
 &\qquad + K^\mathrm{E}\langle \bw_h^{\mathrm{E},\ell+1},\bz_h^\mathrm{E}\rangle_\Sigma - G(\bw_h^{\mathrm{E},\ell+1};\bz_h^\mathrm{E})  = \langle\cR_h^{\mathrm{D},\ell+1},\bz_h^\mathrm{E}\rangle_\Sigma   \quad \forall \bz_h^\mathrm{E}\in\bZ_h^\mathrm{E},
\end{cases}\label{galerkin-wE}\\[2ex]
&\begin{cases}
 a(\bu_h^{\mathrm{E},\ell+1},\bv_h^\mathrm{E})   + b(p_h^{\mathrm{E},\ell+1},\bv_h^\mathrm{E}) + L^\mathrm{E}\langle \bu_h^{\mathrm{E},\ell+1},\bv_h^\mathrm{E}\rangle_{\Sigma\cup \Gamma} \!\!\!\! & = F(\bw_h^{\mathrm{E},\ell+1}; \bv_h^\mathrm{E})+ \langle\cS_h^{\mathrm{D},\ell+1},\bv_h^\mathrm{E}\rangle_\Sigma \ \forall \bv_h^\mathrm{E}\in \bV_h^\mathrm{E},\\
b\bigl(q_h^\mathrm{E},\bu_h^{\mathrm{E},\ell+1}\bigl) -\, c(p_h^{\mathrm{E},\ell+1},q_h^\mathrm{E})  & =  0 \quad \forall q_h^\mathrm{E}\in \cQ_h^\mathrm{E}.
\end{cases}\label{galerkin-upE}
\end{align}

Note that the order in which these blocks are solved is arbitrary. We also note 
that the blocks \eqref{galerkin-wD},\eqref{galerkin-upE} still constitute nonlinear problems, that we solve through inner Newton iterations. Alternatively, one can commute the domain decomposition method and the nonlinear solvers, in such a way that the nonlinearly coupled set of equations can be solved monolithically on each subdomain. Further details on the choice we adopt here, are provided in Algorithm \ref{alg:masteralg}, below. Finally we remark that the sequential steps \eqref{galerkin-wD}-\eqref{galerkin-upE} can be regarded as 
a block Gauss-Seidel iteration applied to the fully monolithic differential-algebraic system.

For the discretisation of the elasticity equations in mixed 
form, we require the discrete spaces in \eqref{eq:subspaces} to be inf-sup stable. In our implementation we 
choose the so-called MINI elements \cite{arnold84}, characterised by 
$$ \bV^*_h := \{\bv \in \bV^*: \bv|_K \in [\PP_1(K) + \mathbb{B}(K)]^d,\ K\in\cT_h^*\} , \quad 
\cQ_h^* := \{ q\in \cQ^*: q|_K \in \PP_1(K),\ K\in \cT_h^*\},  $$ 
where $\mathbb{B}(K)$ is the space of cubic bubble functions defined locally on an element $K$; 
and we employ the space of conforming Lagrangian finite elements for the approximation of the 
species concentrations $\bZ_h^* :=\{ \bz \in\bZ^*: \bz|_K\in\PP_1(K)^m,\ K\in\cT_h^*\}.$

The discretisation \eqref{galerkin-wD}-\eqref{galerkin-upE} leads to the following 
system of equations
\begin{align}
\bsM^* \dot{\mathbf{w}}^{*,l+1} + \bsCu^* \vstep{w}{}{*,l+1} + \bsD \vstep{w}{}{*,l+1} - \bG^{*}(\vstep{w}{}{*,l+1}) &= \cero,\label{gal-Mat-w}\\
\begin{cases}
\bsA^* \vstep{u}{}{*,l+1} + (\bsB^*)^T \vstep{p}{}{*,l+1} - \bF^{*}(\vstep{w}{}{*,\ddagger}) &= \cero, \\
\bsB^* \vstep{u}{}{*,l+1} - \bsC \vstep{p}{}{*,l+1} &= \cero,
\end{cases}\label{gal-Mat-up}
\end{align}
where the superscript $\ddagger$ may be replaced by the appropriate index associated 
to the Schwarz algorithm. We readily observe that only \eqref{gal-Mat-w} is a nonlinear 
ODE system, whereas the matrix form of the elasticity sub-problem consists of 
a linear system. The subscript $u$ in $\bsCu^*$ emphasises that 
the matrix is directly dependent on displacements.

\subsection{Adaptive timestep scheme}\label{sec:adaptive}
The family of fully-discrete methods that we use here is based on adaptive, embedded Runge-Kutta schemes of order 2 and 3 (RK23). For a generic ODE system $\frac{d}{dt} v_h = L_h (v_h, \varphi_h(t))$, we recall that the main step in the classical RK23 scheme reads 
\begin{align*}
v_h^{(i)} &= v_h^n + \sum_{j=1}^{s} a_{ij}\Delta t^n L_h(v_h^{(j)},
\varphi_h(t^n + c_j \Delta t^n)) \quad i = 1,2,3, \\
v_h^{n+1} &= v_h^n + \sum_{j=1}^{s} b_{j}\Delta t^n L_h(v_h^{(j)},
\varphi_h(t^n + c_j \Delta t^n)), \\
\widehat{v}_h^{n+1} &= \widehat{v}_h^n + \sum_{j=1}^{s} \widehat{b}_{j}\Delta t^n L_h(v_h^{(j)},
\varphi_h(t^n + c_j \Delta t^n)),
\end{align*}
with $s=3$ stages, and where the two approximations $v_h,\widehat{v}_h$ are 
of order $\mathcal{O}(p)$ and $\mathcal{O}(\widehat{p})$, respectively. As in 
classical embedded RK methods, here the coefficients $\{a_{ij},c_j\}$ coincide 
in both approximations; however the coefficients $b_j$ differ in order to 
attain the expected convergence orders (in general one chooses $\widehat{p}=p\pm1$, 
cf. \cite{HairerBookODEI}). Even when the stability properties of the two schemes 
do not match, embedded RK methods provide a rather non-expensive 
estimator for the local error (of the less accurate solution) 
which is then directly used in the control of the timestep 
\cite{HairerBookODEI, HairerBookODEII}. The specific description of the 
controller used here is presented in Algorithm \ref{algo:TRBDF2}, below.

\begin{table}[!h]
\begin{center}
   \begin{tabular}{ c | c c c }
     0 & 0 & 0 & 0 \\ 
     $2\gamma$ & $\gamma$ & $\gamma$ & 0 \\
     1 & $1-b_2-\gamma$ & $b_2$ & $\gamma$ \\ \hline
     & $1-b_2-\gamma$ & $b_2$ & $\gamma$ $\vphantom{\int^X}$\\ \hline
     & $1-\widehat{b}_2-\widehat{b}_3$ & $\widehat{b}_2$ & $\widehat{b}_3$  $\vphantom{\int^{X^X}}$
   \end{tabular}
 \end{center}

\vspace{-3mm}
 \caption{Butcher's tableau for the TR-BDF2 method.}\label{tab:TRBDF2}
\end{table}

The RK23 scheme is used together with the TR-BDF2 method (see \cite{bankTRBDF2}), 
 which is a particular class of Diagonally Implicit Runge-Kutta (DIRK) methods (that is, 
$a_{ij}=0$ for $i<j$, and at least one $a_{ii}\neq 0$) \cite{HairerBookODEI}. The 
scheme is characterised by the following Butcher's tableau reported in Table~\ref{tab:TRBDF2},  
where $b_2 = \frac{1-2\gamma}{4\gamma}$, 
$\widehat{b}_2 = \frac{1}{12\gamma(1-2\gamma)}$, 
$\widehat{b}_3 = \frac{2-6\gamma}{6(1-2\gamma)}$ and $\gamma = \frac{2-\sqrt{2}}{2}$.

As established in \cite{HOSEA1996}, the A-stability, the stiffly accuracy ($a_{sj}=b_j$, for $j=1,\ldots,s$), 
as well as the strong S-stability of the method are guaranteed by the coefficients used in Table~\ref{tab:TRBDF2}. An appealing advantage of this scheme is that, for each step, only two implicit stages are required (where one can reuse the same Newton iteration matrix) while one stage is explicit. In this regard, the method is 
closely related to singly DIRK (SDIRK) methods (where all diagonal $a_{ii}$ coefficients are equal and the first stage is explicit) \cite{HOSEA1996}. In addition, TR-BDF2 falls in the class of First-Same-As-Last (FSAL) RK schemes, i.e. methods where the current explicit stage coincides with the last stage 
of the previous step, therefore decreasing the overall computational cost.

Based on Table \ref{tab:TRBDF2}, and considering 
$L_h = \bsCu^* \vstep{w}{}{*} + \bsD \vstep{w}{}{*} - \bG^{*}(\vstep{w}{}{*})$ (where 
we have removed the Schwarz index whenever clear from the context), we can recast the intermediate stages as 
\begin{align*}
\vstep{w}{}{*,(1)} &= \vstep{w}{}{*,n}, \\
\bsM^* \mathbf{w}^{*,(2)} &= \Delta t^n \gamma \bigl[ -\bsCu^* \vstep{w}{}{*,(2)} - \bsD \vstep{w}{}{*,(2)} + \vstep{G}{}{*}(\vstep{w}{}{*,(2)}) \bigr] \\
&\qquad + \Delta t^n \gamma \bigl[ -\bsCu^* \vstep{w}{}{*,(1)} - \bsD \vstep{w}{}{*,(1)} + \vstep{G}{}{*}(\vstep{w}{}{*,(1)}) \bigr] + \bsM^* \mathbf{w}^{*,n}, \\
\bsM^* \mathbf{w}^{*,(3)} &= \Delta t^n \gamma \bigl[ -\bsCu^* \vstep{w}{}{*,(3)} - \bsD \vstep{w}{}{*,(3)} + \vstep{G}{}{*}(\vstep{w}{}{*,(3)})\bigr] \\
&\qquad +  \Delta t^n b_2 \bigl[ - \bsCu^* \vstep{w}{}{*,(2)} - \bsD \vstep{w}{}{*,(2)} + \vstep{G}{}{*}(\vstep{w}{}{*,(2)})\bigr] \\
&\qquad + \Delta t^n (1-b_2-\gamma) \bigl[ -\bsCu^* \vstep{w}{}{*,(1)} - \bsD \vstep{w}{}{*,(1)} + \vstep{G}{}{*}(\vstep{w}{}{*,(1)})\bigr] + \bsM^* \mathbf{w}^{*,n}. 
\end{align*}
These problems should be solved through Newton iterations before moving to the 
next timestep. Alternatively, and as suggested in \cite{HairerBookODEII}, one can apply the change of 
variables $v_{z,h}^{(i)} := v_h^{(i)} - v_h^n$ and rewrite the problem accordingly. 
Doing so results in decreasing round-off errors and serves in reducing the cost to compute $\widehat{v}_h^{n+1}$. We follow this approach and modify the superscripts 
$(1) \rightarrow n$, $(2) \rightarrow \ng$ and $(3) \rightarrow new$, to obtain the following full discretisation of the ODE system:
\begin{align*}
\vstep{w}{}{*,\ng} &= \gamma \vstep{w}{z}{*,\ng} 
+ \gamma \vstep{w}{z}{*,n} + \mathbf{w}^{*,n},  \\
\mathbf{w}^{*,new} &= \gamma \vstep{w}{z}{*,new} 
+ b_2 \vstep{w}{z}{*,\ng} + (1-b_2-\gamma)\vstep{w}{z}{*,n} + \mathbf{w}^{*,n},
\end{align*}
with 
\begin{align*}
\vstep{w}{z}{*,\ng} &= (\bsM^*)^{-1}\Delta t^n \left[ -\bsCu^* \vstep{w}{}{*,\ng} - \bsD \vstep{w}{}{*,\ng} + \vstep{G}{}{*}(\vstep{w}{}{*,\ng})\right],  \\
\vstep{w}{z}{*,new} &= (\bsM^*)^{-1}\Delta t^n \left[ -\bsCu^* \vstep{w}{}{*,new} - \bsD \vstep{w}{}{*,new} + \vstep{G}{}{*}(\vstep{w}{}{*,new})\right].
\end{align*}
This nonlinear system is solved by means of Newton steps as follows 
\begin{align*}
\vstep{w}{z,k+1}{*} &= \vstep{w}{z,k}{*} + \dvstep{k+1}{}{w}{z}{*}, \\
\biggl(\frac{1}{\dtnum{n}}\bsM^*-\gamma\bsJ^*\biggr)^{-1} \dvstep{k+1}{}{w}{z}{*} &= \left[ -\bsCu^* \vstep{w}{k}{*} - \bsD \vstep{w}{k}{*} + \vstep{G}{}{*}(\vstep{w}{k}{*})\right] - \frac{1}{\Delta t^n}\bsM \vstep{w}{z,k}{*},
\end{align*}
where 
\begin{equation*}
	\bsJ^* = -\bsD^* + \frac{\partial \vstep{G}{}{*}(\vstep{w}{k}{*})}{\partial \vstep{w}{}{*}} - \bsC_{u}^* = -\bsD^* + \bsJ_{\bG^*} - \bsC_{u}^*.
\end{equation*}
We finally point out that an additional reduction of computational 
cost occurs due to the fact that the local error 
$\mathrm{err}=\widehat{v}_h^{n+1} - v_h^{n+1}$ can be accessed by the simple 
relation 
$$\mathrm{err}=(\widehat{b}_1 - b_1)\vstep{w}{z}{n} + (\widehat{b}_2-b_2)\vstep{w}{z}{\ng}+(\widehat{b}_3-b_3)\vstep{w}{z}{new},$$ 
avoiding the explicit storage and computation of $\widehat{v}_h^{n+1}$.

\subsection{Algorithmic details}


\textbf{Overall solution algorithm.} The global numerical scheme to solve the mechanochemical model  
is summarised in \autoref{alg:masteralg}. After assigning suitable initial conditions 
$ \vstep{u}{}{*,0} = \cero$, $\vstep{p}{}{*,0}=\cero$, $\vstep{w}{}{*,0} = \vstep{w}{0}{*}$,  
projected to the corresponding finite element spaces, we use the non-overlapping Schwarz 
method \eqref{galerkin-wD}-\eqref{galerkin-upE} to produce a solution until time $T_f$. 
Two important points stand out from this master algorithm. First, only the ADR 
equations are solved through the embedded RK method within the TR-BDF2 method 
(see \autoref{alg:masteralg}, \autoref{MA:TRBDF2}, detailed in \autoref{algo:TRBDF2}), whereas 
the elasticity subproblem is solved by an implicit method (cf. \autoref{MA:elast} in \autoref{alg:masteralg}). 
Furthermore, the velocity in the advection term is approximated from the solid displacement using a backward Euler scheme (cf. \autoref{alg:masteralg}, \autoref{MA:alg1-vel1} and \ref{MA:alg1-vel2}). Second, the computation performed in \autoref{alg:masteralg}, \autoref{MA:rescale1} and \ref{MA:rescale2}, represents an alternative way of performing the explicit first stage of the RK scheme \cite{HOSEA1996}. Such a rescaling must be executed at each modification of the timestep, so that one starts with an adequate value for $\vstep{w}{z}{}$ when entering \autoref{MA:TRBDF2} of \autoref{alg:masteralg}. An additional benefit of the rescaling is that it improves the stability of the overall stiff problem, and at the same time avoids unnecessary function evaluations.

\begin{algorithm}[t] 
\caption{Overall solution algorithm for the mechanochemical problem} 
\label{alg:masteralg}
{\small\begin{algorithmic}[1]

\State Initialise the solution vectors $\vstep{u}{}{*,0}$, $\vstep{p}{}{*,0}$, $\vstep{w}{}{*,0}$ 
\And \textbf{define} $\vstep{w}{}{} \triangleq \left[ \vstep{w}{}{D}, \vstep{w}{}{E} \right]^T$, $\vstep{w}{z}{} \triangleq \left[ \vstep{w}{z}{D}, \vstep{w}{z}{E} \right]^T$
\While{$t \leq T_f$}
\State $\dtnum{\min} = \epsilon|t|$, $\dtnum{\mathrm{old}} = \dtnum{}$, $\dtnum{} = \min( \dtnum{\max}, \max( \dtnum{\min}, \dtnum{\mathrm{old}} ) )$
	\If{ $\dtnum{}\neq\dtnum{\mathrm{old}}$ }
		\State $\vstep{w}{z}{n} = \frac{\dtnum{}}{\dtnum{\mathrm{old}}} \vstep{w}{z}{n}$\label{MA:rescale1}
		\State \texttt{needNewMJ} $\leftarrow$ \True
	\EndIf
	\State \textbf{Compute} $\tn{\mathrm{new}},\vstep{w}{}{\mathrm{new}},\vstep{w}{z}{\mathrm{new}}$ using \Call{TR-BDF2}{$\vstep{w}{}{n}$,$\vstep{w}{z}{n}$,$\dtnum{}$,$\bbeta_{\dtnum{}}$} \Comment{See Algorithm \ref{algo:TRBDF2}} \label{MA:TRBDF2}
	\For{$*\in\{D,E\}$} \textbf{procedure} \textsc{ElasticitySolver} associated to \eqref{galerkin-upD} and \eqref{galerkin-upE}\label{MA:ForDE}
			\State $\bF^{*,n+1}=\bF(\vstep{w}{}{*,new})$, $\vstep{v}{}{*,n+1}=-\frac{\vstep{u}{}{*,n}}{\dtnum{}}$ \label{MA:alg1-vel1}
			\State $\left[ \begin{array}{c} \vstep{u}{}{*,n+1} \\ \vstep{p}{}{*,n+1} \end{array}\right] = \left[\begin{array}{cc} \bsA^* & (\bsB^*)^T \\ \bsB^* & -\bsC^* \\ \end{array}\right]^{-1}\left[ \begin{array}{c} \bF^{*,n+1}\\ \cero \end{array}\right]$ \label{MA:elast}
			\State $\vstep{v}{}{*,n+1}=\vstep{v}{}{*,n+1} + \frac{\vstep{u}{}{*,n+1}}{\dtnum{}}$ \label{MA:alg1-vel2}
			\State \textbf{update} Robin transmission conditions $\cS_h^*$ from \eqref{eq:interface-elast}
	\EndFor
	\State \textbf{set} $n=n+1$, $\tn{n+1} = \tn{\mathrm{new}}$, $\vstep{w}{}{n+1} = \vstep{w}{}{\mathrm{new}}$, $\vstep{w}{z}{n+1}=\vstep{w}{z}{\mathrm{new}}$, \And 
		 $| \vstep{w}{}{n+1} |_2 = | \vstep{w}{}{\mathrm{new}} |_2$
		\State \textbf{Compute} $q=\left(\frac{\mathrm{err}}{\mathrm{R_{TOL}}}\right)^{k_I}$, $\mathrm{\mathrm{ratio}}=\frac{\dtnum{\max}}{\dtnum{}}$, 
		 $\mathrm{ratio}=\min(\mathrm{ratiomax},\max(\mathrm{ratiomin},\mathrm{ratio}))$
		\If{$|\mathrm{ratio}-1|>\mathrm{ratiomin}$}
			\State $\dtnum{}=\mathrm{ratio}\cdot\dtnum{}$, $\vstep{w}{z}{n+1}=\mathrm{ratio}\cdot\vstep{w}{z}{n+1}$\label{MA:rescale2}
			\State \texttt{needNewMJ} $\leftarrow$ \True
		\EndIf
\EndWhile
\end{algorithmic}
}
\end{algorithm}


\textbf{Adaptive Runge-Kutta method.} The RK procedure is detailed in \autoref{algo:TRBDF2} and it 
features a similar structure as the ODE solver \textsf{ode23tb} in MATLAB. Again, 
we highlight two major components: the one-step RK resolution and the adaptive setting of the timesteps. \autoref{TRBDF2:S1} in \autoref{algo:TRBDF2} represents stage 2 of the RK method, i.e., the first implicit stage (\texttt{S1}). As the scheme is implicit, we initialise the tangent system associated to the Newton iteration. These matrices will actually be assembled only if the Boolean \texttt{mjcontrol} is true, and this switch is described in  \autoref{algo:mjcontrol}, below. 
Should \texttt{mjcontrol} be false, the mass-Jacobian blocks are evaluated at each Newton step (see \autoref{algo:TRBDF2}, specifically \autoref{TRBDF2:RHSS1}). 
We also stress that the initial guess for the rescaled derivative $\vstep{w}{z}{\ng}$ can be simply 
taken as $\vstep{w}{z}{n}$, since extrapolating the derivative of the interpolant from the previous step does not improve accuracy nor efficiency. Next, when updating the solutions we also change the Robin transmission conditions so that the Schwarz algorithm is performed at the same time (see \autoref{TRBDF2:RobS1} in \autoref{algo:TRBDF2}). 
In order to retrieve $\vstep{w}{z,k+1}{new}$ and $\vstep{w}{k+1}{new}$, a second Newton step is necessary within 
the implicit stage  (\texttt{S2}), the particularity of which is the evaluation of the initial guess $\vstep{w}{z,0}{new}$, by extrapolation of the derivatives of the cubic Hermite interpolant from $t^n$ to $t^{\ng}$ \cite{HOSEA1996}. 

\begin{algorithm}[!t] 
\caption{Adaptive Runge-Kutta time stepping: 
\textbf{procedure} \textsc{TR-BDF2}($\vstep{w}{}{n}$,$\vstep{w}{z}{n}$,$\dtnum{}$,$\bbeta_{\dtnum{}}$)} 
\label{algo:TRBDF2}
{\small
\begin{algorithmic}[1] 
\State \textbf{compute} $\mathrm{sc}_i = \max( | \vstep{w}{}{n} |_2, \eta )$\label{TRBDF2:sc}
		\State \textbf{assemble} $\bM$ and $\bJ$ using \Call{MJcontrol}{$\vstep{w}{}{n}$,$\vstep{w}{z}{n}$,$\dtnum{}$,$\bbeta_{\dtnum{}}$} \Comment{See Algorithm \ref{algo:mjcontrol}}\label{TRBDF2:S1}
		\State \textbf{set} $\tn{\ng} = t + 2\gamma\dtnum{}$, $\vstep{w}{0}{\ng} = \vstep{w}{}{n} + 2\gamma\vstep{w}{z}{n}$, $\vstep{w}{z,0}{\ng} = \vstep{w}{z}{n}$ \Comment{\textbf{Stage 2}}
		\State \texttt{rejectS1} $\leftarrow$ \False
		\State \textbf{set} $\res{S_2,min} = 100\epsilon\vvvert\vstep{w}{0}{\ng}\vvvert$
		\For{$k = 0, \ldots, K$}\label{TRBDF2:NewtS1}
			\For{$*\in\{D,E\}$}\label{TRBDF2:forDE}
				\State {$\vstep{g}{k}{*,\ng} = -\frac{1}{\Delta t}\bsM^*\vstep{w}{z,k}{*,\ng} + \left[ -\bsD^*\vstep{w}{k}{*,\ng}+\bG^*(\vstep{w}{k}{*,\ng})-\bsCu^{*,n}\vstep{w}{k}{*,\ng} \right]$}\label{TRBDF2:RHSS1}
				\If{$\neg$ \texttt{mjcontrol}}
				\State {Update $\bsM^*=\bsM_k^{*,\ng}$ \And $\bJ^*=\bJ_k^{*,\ng}$}
				\EndIf
				\State \textbf{solve} {$\dvstep{k+1}{\ng}{w}{z}{*} = (\MJs)^{-1}\vstep{g}{k}{*,\ng}$}
				\State \textbf{compute} $\vstep{w}{z,k+1}{*,\ng} = \vstep{w}{z,k}{*,\ng} + \dvstep{k+1}{\ng}{w}{z}{*}$, $\vstep{w}{k+1}{*,\ng} = \vstep{w}{k}{*,\ng} + \gamma \dvstep{k+1}{\ng}{w}{z}{*}$		
				\State \textbf{update} Robin transmission conditions $\cR_h^*$ from \eqref{eq:interface-adr}\label{TRBDF2:RobS1}
			\EndFor	
			\State \Call{NewtonConv}{$\vstep{w}{k+1}{\ng}$,$\dvstep{k+1}{\ng}{w}{z}{}$,$\res{S_2,min}$,$\bbeta_{\dtnum{}}$} \Comment{See  \autoref{algo:NewtonConv}}\label{TRBDF2:NConvS1}
		\EndFor
		\State \texttt{rejectS2} $\leftarrow$ \False
		\If{ $\neg$ \texttt{rejectS1}}\Comment{\textbf{Stage 3}}\label{TRBDF2:S2}
				\State \textbf{compute} $\mathrm{sc}_i = \max( | \vstep{w}{}{\ng} |_2, \mathrm{sc}_i )$
			\State $\vstep{w}{z,0}{\mathrm{new}} = p_{31}\vstep{w}{z}{n} + p_{32}\vstep{w}{z}{\ng} + p_{33}\vstep{w}{}{\ng} - p_{33}\vstep{w}{}{n}$
			\State $\vstep{w}{0}{\mathrm{new}} = \vstep{w}{}{n} + (1-b_2-\gamma)\vstep{w}{z}{n} + b_2\vstep{w}{z}{\ng}+\gamma\vstep{w}{z,0}{\mathrm{new}}$
			\State \textbf{redo} \underline{steps 6-15} to retrieve $\vstep{w}{z,k+1}{\mathrm{new}}$ \And $\vstep{w}{k+1}{\mathrm{new}}$  
		\EndIf
		\If{\texttt{rejectS1} $\lor$ \texttt{rejectS2}} \label{TRBDF2:RejStage}
			\If{\texttt{mjcontrol}}
				\If{\texttt{Jcurrent} $\land$ \texttt{Mcurrent}}
						\State \textbf{set} $\dtnum{\mathrm{old}}=\dtnum{}$, $\dtnum{}=\max(fac_{S1S2}\dtnum{\mathrm{old}},\dtnum{\min})$, 
						 $\vstep{w}{z}{n}=\frac{\dtnum{}}{\dtnum{\mathrm{old}}}\vstep{w}{z}{n}$
						\State \texttt{needNewMJ} $\leftarrow$ \True
				\Else
					\State \texttt{needNewJ} $\leftarrow$ $\neg$\texttt{Jcurrent}, \texttt{needNewM} $\leftarrow$ $\neg$\texttt{Mcurrent}\label{TRBDF2:MJCur}
				\EndIf
			\Else
					\State \textbf{set} $\dtnum{\mathrm{old}}=\dtnum{}$, $\dtnum{}=\max(fac_{S1S2}\dtnum{\mathrm{old}},\dtnum{\min})$, 
					 $\vstep{w}{z}{n}=\frac{\dtnum{}}{\dtnum{\mathrm{old}}}\vstep{w}{z}{n}$
			\EndIf
		\Else\label{TRBDF2:RejStep}
		\State \textbf{compute} $\mathrm{sc}_i = \max( | \vstep{w}{}{\mathrm{new}} |_2, \mathrm{sc}_i )$
			\State $\ewo = (\hat{b}_1 - b_1)\vstep{w}{z}{n} + (\hat{b}_2 - b_2)\vstep{w}{z}{\ng} + ( \hat{b}_3 - b_3)\vstep{w}{z}{n+1}$
			\If{$\neg$ \texttt{mjcontrol}}
				\State Update $\bM=\bM^{\mathrm{new}}$ \And $\bJ=\bJ^{\mathrm{new}}$
			\EndIf
			\State \textbf{compute} { $\ewt = (\MJ)^{-1}\ewo$ } \And  
			 $\mathrm{err}=\max(\vvvert\ewt\vvvert,\frac{\vvvert\ewo\vvvert}{16})$\label{TRBDF2:e2w}
			\If{$\mathrm{err}>\mathrm{R_{TOL}}$}
				\State $\dtnum{\mathrm{old}}=\dtnum{}$, $\dtnum{}=\max\left(\dtnum{\mathrm{old}}\max\left(fac_{\min},fac\left(\frac{\mathrm{R_{TOL}}}{err}\right)^{k_I}),\dtnum{\min}\right)\right)$, 
				 $\vstep{w}{z}{n}=\frac{\dtnum{}}{\dtnum{\mathrm{old}}}\vstep{w}{z}{n}$
				\State \texttt{needNewMJ} $\leftarrow$ \True
			\EndIf
		\EndIf
\end{algorithmic}}
\end{algorithm}

The algorithm also focuses on the automatic determination of each timestep, 
in view of the stability of the time integration of stiff problems. The process 
is designed to reject the overall step if one of the two implicit stages fails or if the local error is too large. The potential rejection at an implicit stage is monitored by the Boolean variables \texttt{rejectS1} and \texttt{rejectS2}, which are switched on according to the value of the residual in the Newton iteration 
described in \autoref{algo:NewtonConv}. How the value is actually computed 
may affect the performance of the algorithm (see e.g. \cite{Söderlind2002}). 
For a generic residual vector $\bv$ of length $M$, here we adopt the rescaled maximum norm 
$$\vvvert \bv \vvvert = \max_{i=1,\ldots,M}(v_i\cdot sc^{-1}_i),$$
where $sc_i:= \max( sc_i, |w_{k+1,i}|)$ is the local maximum of the solution 
vector at the $k$-th Newton step, and it is reinitialised as 
$sc_i = \max( |\vstep{w}{}{n}|_2,\eta)$, where $\eta = \frac{\mathrm{A_{TOL}}}{\mathrm{R_{TOL}}}$. 
This process is summarised in \autoref{TRBDF2:sc} from \autoref{algo:TRBDF2}. 

We recall (see \autoref{NC:res} of \autoref{algo:NewtonConv}) that 
the residual in the Newton step is simply the norm of the current  
increment $\dvstep{k+1}{}{w}{z}{*}$. The loop can break under different 
cases. For instance, if the residual is smaller than a given threshold, 
then \autoref{NC:Conv1} terminates the step. Otherwise, the linear 
convergence of the overall scheme (see e.g. \cite{HairerBookODEII}) implies that 
one can evaluate the iteration error using the residual together with the 
inequality $\|\dvstep{k+1}{}{w}{z}{}\|\leq\theta\|\dvstep{k}{}{w}{z}{}\|$, with $\theta<1$.
This experimental rate is reset to zero at each mass-stiffness assembly so that 
it depends on the Newton iteration matrix \cite{remani2012}. 
The rate is then estimated using $\frac{\theta}{1-\theta}\mathrm{res}$ (see \autoref{NC:TE}) with two successive iterations and the current error. We consider that the algorithm has converged if the error is smaller than 
 $\kappa_N$ times a fixed tolerance. 
 A given stage is therefore rejected if at least one of the following situations occurs: 
(a) the convergence rate is larger than 1 (cf. \autoref{NC:Div1}); (b)
the maximum number of iterations $K=10$ is reached (cf. \cite{HairerBookODEII} and \autoref{NC:DivIt}); and 
(c) For any $k$, we have $\kappa_N\mathrm{TOL}_N < \eN \theta^{K-k}$ (cf. \autoref{NC:Div2}).
Point (c) represents a rough estimate of the iteration error expected after $K-1$ iterations.


\begin{algorithm}[t] 
\caption{Newton stopping criterion: \textbf{procedure} \textsc{NewtonConv}($\vstep{w}{k+1}{}$,$\dvstep{k+1}{}{w}{z}{}$,$\res{\min}$,$\bbeta_{\dtnum{}}$)} 
\label{algo:NewtonConv}
{\small 
\begin{algorithmic}[1]
\State \textbf{compute} $\widehat{\mathrm{sc}}_i = \max( \mathrm{sc}_i, |\vstep{w}{k+1,i}{}| )$\label{NC:res}
	\And $\res{k+1}(\widehat{\mathrm{sc}}_i) = \vvvert \dvstep{k+1}{}{w}{z}{} \vvvert$
		\If{$\res{k+1} \leq \res{\min}$}\label{NC:Conv1}
			\State \Break
		\ElsIf{$k=0$}\label{NC:k0}
			\If{\texttt{needNewRate} $\land$ Stage 2}\label{NC:NewRate}
				{\If{\texttt{mjcontrol}}
					\State \texttt{needNewRate} $\leftarrow$ \False
				\EndIf}
				\State $\theta=0$
			\Else\label{NC:OldRate}
				\State $\eN = \frac{\theta}{1-\theta}\res{k+1}$
				\If{$\eN \leq 0.1\kappa_N \mathrm{TOL}_N$}
					\State \Break
				\EndIf
			\EndIf
		\ElsIf{$\res{k+1}>0.9\res{k}$}\Comment{Divergence}\label{NC:Div1}
			\State \texttt{rejectSX} $\leftarrow$ \True\Comment{\texttt{X} stands for the stage}
			\State \Break
		\Else
			\State \textbf{set} $\theta = \max\left( 0.9\theta, \frac{\res{k+1}}{\res{k}} \right)$
			\And $\eN = \frac{\theta}{1-\theta}\res{k+1}$\label{NC:TE}
			\If{$\eN\leq\kappa_N\mathrm{TOL}_N$}\Comment{Convergence}\label{NC:Conv2}
				\State \Break
			\ElsIf{$k=K$} \Comment{Max. iterations attained}\label{NC:DivIt}
				\State \texttt{rejectSX} $\leftarrow$ \True
				\State \Break
			\ElsIf{$\kappa_N\mathrm{TOL}_N < \eN \theta^{K-k}$}\Comment{Divergence}\label{NC:Div2}
				\State \texttt{rejectSX} $\leftarrow$ \True
				\State \Break
			\EndIf
		\EndIf
\end{algorithmic}}
\end{algorithm}

Should a stage be rejected, one moves to  \autoref{TRBDF2:RejStage} in \autoref{algo:TRBDF2} 
and reduces the timestep $\Delta t$ by a factor $fac_{S1S2}$ (independently of the 
mass matrix assembly). Then, \autoref{TRBDF2:MJCur} checks whether both 
mass and Jacobian matrices have been updated. Finally, \autoref{TRBDF2:RejStep} verifies 
that the local error estimation remains below a prescribed tolerance. 
Based on the estimators for the scaled derivatives we proceed to compute the local asymptotic error 
$\vstep{e}{w}{1}$ as a result from the linear system stated in \autoref{TRBDF2:e2w} of \autoref{algo:TRBDF2}. 
This step is necessary as the boundedness of the raw estimator is not always guaranteed \cite{HairerBookODEII}. 

Next, we select the maximum error between these two estimators as an approximation of the timestep error. 
In contrast with \cite{HairerBookODEII}, we do not observe here the so-called \emph{hump} phenomenon 
intrinsic to automatic step methods, so we do not require an additional correction step. We only update 
the timestep in an automatic fashion (thanks to the sequence of rejection procedures), capturing correctly 
the dynamics of the integrated system. When all errors are lower than the corresponding tolerance, the 
timestep is accepted and the interface elasticity problem is solved in \autoref{MA:elast} from 
\autoref{alg:masteralg}. We stress that the algorithm will stop if $\dtnum{}\leq\dtnum{\min}$.


\begin{algorithm}[!t] 
\caption{Control of mass-stiffness assembly: \textbf{procedure} \textsc{MJcontrol}($\vstep{w}{}{n}$,$\vstep{w}{z}{n}$,$\dtnum{}$,$\bbeta_{\dtnum{}}$) } 
\label{algo:mjcontrol}
{\small\begin{algorithmic}[1]
		\If{\texttt{needNewM}$\land$\texttt{mjcontrol}}
			\State {\textbf{compute} $\bsM^{*,n}$, with $*\in\{D,E\}$}
			\State \texttt{Mcurrent}$\leftarrow$\True, \texttt{needNewM}$\leftarrow$\False, \texttt{needNewMJ}$\leftarrow$\True 
		\EndIf
		\If{\texttt{needNewJ}$\land$\texttt{mjcontrol}}
			\State \textbf{compute} {$\bsJ^{*,n}=-\bsD^*+\bsJ_{\bG^*}^n-\bsC_{u}^{*,n}$ with $*\in\{D,E\}$}
			\State \texttt{Jcurrent}$\leftarrow$\True, \texttt{needNewJ}$\leftarrow$\False, \texttt{needNewMJ}$\leftarrow$\True 
		\EndIf
		\If{\texttt{needNewMJ}$\land$\texttt{mjcontrol}}
			\State \textbf{compute} {$\frac{1}{\Delta t}\bsM^{*,n}-\gamma\bsJ^{*,n}$ with $*\in\{D,E\}$}
			\State \texttt{needNewMJ}$\leftarrow$\False, \texttt{needNewRate}$\leftarrow$\True
		\EndIf
\end{algorithmic}}
\end{algorithm}

\autoref{algo:mjcontrol} permits to control the mass-stiffness assembly. This procedure consists of a 
number of controllers, including the aforementioned \texttt{mjcontrol} (see a list in Table~\ref{tab:boolean}). 
The state modification of these variables depends on step-size revaluation, which itself depends on the 
rejected stage of the RK 
method. Notice that  \texttt{Mcurrent} and \texttt{Jcurrent} control whether the mass and Jacobian matrices 
are evaluated at each new timestep (if this timestep is rejected at least once). For models such as \eqref{eq:rd},
the label \texttt{Mcurrent} is always \True. However if we want to solve the interface problem 
using e.g. a level-set approach, then the mass matrix needs to be recomputed and the controller can 
be exploited.

\begin{table}[!h]
\renewcommand{\arraystretch}{1.25}
\begin{center}
   \begin{tabular}{ ||K{0.2\textwidth} | K{0.7\textwidth} ||}
   \hline
\hline
     \texttt{mjcontrol} & Switch to reuse Newton iteration matrices over two implicit stages \\ 
     \texttt{needNewM}  & Need to evaluate mass matrix $\bM^*$ \\
     \texttt{needNewJ}  & Need to evaluate Jacobian matrix $\bJ^*$ \\
     \texttt{needNewMJ} & Need to evaluate mass-Jacobian matrix $\MJs$  \\
     \texttt{Mcurrent}  &Appropriate mass matrix $\bM$ for the current timestep \\
     \texttt{Jcurrent}  & Appropriate Jacobian matrix $\bJ$ for the current timestep\\
\hline
\hline
   \end{tabular}
 \end{center}
 
 \vspace{-3mm}
 \caption{Set of Booleans controlling the evaluation and assembly of mass and stiffness matrices.}\label{tab:boolean}
\end{table}

\section{Numerical tests}\label{sec:results}
This section contains a set of examples assessing the experimental 
spatio-temporal convergence  of our primal-mixed method, 
and also illustrates its use in two mechanochemical interface problems. The 
implementation is carried out with the open-source finite element 
library  FreeFem++ \cite{freefem}, and all sparse linear systems are solved with distributed 
direct methods (MUMPS and UMFPACK).

\begin{table}[!h]
\begin{center}
\begin{tabular}{l}
\hline
\hline\noalign{\smallskip}
$c_f^\mathrm{D} = 150$, $c_g^\mathrm{D} = 1$, 
$\rho_0^\mathrm{D} = \rho_2^\mathrm{D}= \rho_3^\mathrm{D} = \rho_5^\mathrm{D}= 1$, 
$\rho_1^\mathrm{D} = 0$, $\rho_4^\mathrm{D} = 0.35$, $M_{11}^\mathrm{D}=1$, $M_{22}^\mathrm{D}=30$,\\[1.5ex] 
$c_f^\mathrm{E} = 20$, $c_g^\mathrm{E} = 2$,  $\rho_0^\mathrm{E} = 2$, 
$\rho_2^\mathrm{E}=\rho_3^\mathrm{E} =2$, $\rho_5^\mathrm{E}= 1$, 
$\rho_1^\mathrm{E} = 0$, $\rho_4^\mathrm{E} = 0.15$, $M_{11}^\mathrm{E}=2$, $M_{22}^\mathrm{E}=10$,\\[1.5ex]
$E^\mathrm{D} = 1000$, $\nu^\mathrm{D}=0.475$, $J^\mathrm{D}=1$, $K^\mathrm{D}=1$e5, 
$E^\mathrm{E} = 10$, $\nu^\mathrm{D}=0.33$, $J^\mathrm{E}=1$, $K^\mathrm{E}=1$e5, 
$\alpha^\mathrm{E}=2.5$ \\
\noalign{\smallskip}\hline
\hline
\end{tabular}
\caption{Example 1. Parameters for the mechanochemical 
model using \eqref{reaction:GM} and \eqref{2way-coupling}.}\label{table:ex01params}
\end{center}
\end{table}

\paragraph{Example 1: Experimental accuracy against manufactured solutions.}  Let us 
consider a rectangular domain $\Om = (0,1)\times(0,1.4)$ divided into the dermis $\OmD = (0,1)^2$ 
and epidermis $\OmE = \Om\setminus\OmD$ subdomains. The initial position of the 
interface is characterised by the segment $x\in (0,1),y=1$, and the external surface 
is $x\in (0,1),y=1.4$. The main goal of our first series 
of tests is to analyse the convergence of the spatial discretisation. We therefore 
focus on a stationary problem where no advection occurs in the morphogen model. 
 We construct smooth exact solutions, defined as follows: 
 \begin{equation}\label{eq:exact-h}
\tilde\bw = \begin{pmatrix} 
 1 - \cos(2\pi x)\sin(3 \pi y) \\
 1 + \frac{1}{2}\cos(2\pi x)\sin(3\pi y)\end{pmatrix}, \quad 
 \tilde\bu = \begin{pmatrix} 
x(1-x)\cos(\pi x) \sin(2\pi y)\\
\sin(\pi x)\cos(\pi y) y^2(1-y)\end{pmatrix},\end{equation}
and set $\buD = \tilde\bu|_{\OmD}$, $\buE = \tilde\bu|_{\OmE}$, $\bw^\mathrm{D} = \tilde\bw|_{\OmD}$, 
$\bw^\mathrm{E} = \tilde\bw|_{\OmE}$, and $p^* = -\lambda^* \vdiv\tilde \bu^*$. 
The functions in \eqref{eq:exact-h} are used to construct the Dirichlet {datum} for the 
displacements on $\partial\Om\setminus\Gamma$ as well as the flux conditions 
imposed on the morphogens at the boundary $\partial\Om$. Suitable non-homogeneous forcing terms are also set using these closed-form solutions. We choose the two-species Gierer-Meinhardt reaction model \eqref{reaction:GM}, and adopt the coupling terms as in \eqref{2way-coupling}. Model parameters are chosen as in Table~\ref{table:ex01params}.

\begin{table}[!h]
\setlength{\tabcolsep}{5pt}
\begin{center}
{\small 
\begin{tabular}{ccccccccccccc}
\hline
\hline
\multicolumn{11}{c}{Space convergence $\vphantom{\int^X}$} \\
\hline\noalign{\smallskip}
DoF & $h$  & $e_1(\bw^\mathrm{D})$ & $r_1(\bw^\mathrm{D})$ & $e_1(\bw^\mathrm{E})$ & $r_1(\bw^\mathrm{E})$ &  $e_0(p^\mathrm{D})$ & $r_0(p^\mathrm{D})$ & $e_0(p^\mathrm{E})$ & $r_0(p^\mathrm{E})$ \\
\noalign{\smallskip}
\hline\noalign{\smallskip}
48 & 0.5207 &  5.6275 & -- & 3.5698 & -- &  1175.8 & --      & 2.2114 & --   \\
108 & 0.2828 &  4.1858 & 0.5794  & 2.6495 & 0.4886 &444.12 & 1.9061 & 1.2433 & 0.9437\\
280 & 0.1736 &  2.6223 & 0.7956  & 1.8701 & 0.7133 &129.18 & 2.1008 & 0.6340 & 1.3790\\
900 & 0.0889 &  1.4380 & 0.9447 & 0.9974 & 0.9398 & 37.750 & 1.9344 & 0.2178 & 1.5973\\
3264 & 0.0432 &  0.7485 & 0.9844  & 0.4852 & 0.9978 & 11.335 & 1.8138 & 0.0801 & 1.3849\\
12276 & 0.0218 &  0.3829 & 0.9890 & 0.2480 & 0.9788 & 3.6282 & 1.6805 & 0.0409 & 0.9624\\
47320 & 0.0110 &  0.1970 & 0.9695 & 0.1326 & 0.9217 & 1.4185 & 1.3701 & 0.0252 & 0.9748 \\
186276 & 0.0055 &  0.1065 & 0.9928 & 0.0684 & 0.9605 &  0.7001 & 0.9599 & 0.0141 & 0.9635\\
\hline\noalign{\smallskip}
DoF & $h$  & $e_0(\buD)$ & $r_0(\buD)$ & $e_1(\buD)$ & $r_1(\buD)$ & $e_0(\buE)$ & $r_0(\buE)$& $e_1(\buE)$ & $r_1(\buE)$ & Iter \\
\noalign{\smallskip}
\hline\noalign{\smallskip}
48 & 0.5207 & 0.0605 & --     & 1.2176 & --     & 0.0423 & --     & 0.5038 & --     & 4 \\
108 & 0.2828 &0.0255 & 1.6925 & 0.6149 & 1.3372 & 0.0192 & 1.2900 & 0.3408 & 0.6409 & 4\\
280 & 0.1736 &0.0081 & 1.9584 & 0.2545 & 1.5009 & 0.0093 & 1.4921 & 0.2226 & 0.8723 & 4\\
900 & 0.0889 &0.0024 & 1.9235 & 0.1201 & 1.1808 & 0.0029 & 1.7335 & 0.1131 & 1.0121 & 4\\
3264 & 0.0432 &0.0006 & 1.9717 & 0.0596 & 1.0563 & 0.0011 & 1.2942 & 0.0537 & 1.0313 & 4\\
12276 & 0.0218 &0.0002 & 1.9373 & 0.0299 & 1.0180 & 0.0006 & 1.9020 & 0.0275 & 0.9765 & 4\\
47320 & 0.0110 &3.61e-5 & 1.9734 & 0.0150 & 1.0061 & 3.70e-5 & 1.9258 & 0.0151 & 0.9848 & 4\\
186276 & 0.0055 &9.02e-6 & 1.8221 & 0.0075 & 1.0020 & 9.12e-6 & 1.8300 & 0.0065 & 0.9730 & 4\\
\hline
\hline
\end{tabular}
\begin{tabular}{cccccc}
\multicolumn{6}{c}{Time convergence $\vphantom{\int^X}$} \\
\hline\noalign{\smallskip}
$\Delta t$  &  $e_{\Delta t}(\bw^\mathrm{D})$ & $r_{\Delta t}(\bw^\mathrm{D})$ & $e_{\Delta t}(\bw^\mathrm{E})$ & $r_{\Delta t}(\bw^\mathrm{E})$ & avg(iter) 
$\vphantom{\int^X}$  \\
\noalign{\smallskip}
\hline
0.0100 & 0.2002 & --     & 0.0881 & -- & 3 \\
0.0050 & 0.0982 & 1.1294 & 0.0292 & 1.1752 & 3 \\
0.0025 & 0.0230 & 1.9301 & 0.0094 & 1.7388 & 3 \\
0.0013 & 0.0056 & 1.9614 & 0.0024 & 1.7929 & 3 \\
0.0006 & 0.0016 & 1.9811 & 0.0006 & 1.9952 & 3 \\
0.0002 & 0.0004 & 1.9937 & 0.0002 & 1.9874 & 3 \\
\noalign{\smallskip}
\hline
\end{tabular}}

\caption{Example 1. Spatial and temporal error history associated to the discretisation 
of the model problem. Errors on the bottom sub-table were obtained on a fine mesh of size $h=0.0013$, and 
computed until a final time $T=0.02$.}\label{table:errors}
\end{center}
\end{table}

We proceed to generate a sequence of successively refined unstructured triangular 
meshes for both subdomains, and compute errors between exact and approximate 
numerical solutions (recall that we are solving for the steady version of the 
coupled problem) on each refinement level. For a generic scalar or vector 
field $s^*$ we will denote errors and convergence rates as follows 
$$e_0(s^*):=\| s^* - s^*_h\|_{0,\Om^*}, \quad e_1(s^*):=\| s^* - s^*_h\|_{1,\Om^*}, 
\quad r_k(s^*)=\frac{\log(e_k(s^*)/\widehat{e_k}(s^*))}{\log(h/\widehat{h})}, \quad k=0,1,$$ 
where $e_k$ and $\widehat{e_k}$ denote errors computed 
on two consecutive meshes of sizes $h$ and $\widehat{h}$, respectively. 
The outcome of the convergence test (using the MINI-element for the approximation 
of displacement and pressure, and Lagrangian finite elements for the species 
concentration) is presented in the first part of Table~\ref{table:errors}, which indicates that 
the method converges optimally. Four Newton iterations are sufficient for each mesh to reach the tolerance of $\mathrm{TOL}_N=1$e-5, and the stopping criterion is based on the $\ell^\infty-$norm of the total residual. 
The time convergence is assessed by taking a fine mesh of size $h=0.0013$ and computing 
approximate solutions to the transient coupled problem \eqref{eq:elasticity},\eqref{eq:rd} 
compared against the exact displacement and pressure from \eqref{eq:exact-h} whereas 
the species concentrations are now $$\bw = \tilde\bw \exp(-kt),$$
with $k= \log(1/2)$. Errors and convergence rates are measured and denoted as follows 
$$e_{\Delta t}(s^*):= \sum_{n = 1}^N \| s^*(t^n) - s^*_h(t^n)\|_{0,\Om^*}, \quad 
r_{\Delta t}(s^*)=\frac{\log(e_{\Delta t}(s^*)/\widehat{e_{\Delta t}}(s^*))}{\log(\Delta t/\widehat{\Delta t})},$$ 
where $e_{\Delta t},\widehat{e_{\Delta t}}$ denote errors computed using two 
consecutive timesteps $\Delta t$ and $\widehat{\Delta t}$, respectively. The 
second block in Table~\ref{table:errors} shows an asymptotic second order convergence, 
as expected from the use of the TR-BDF2 time stepping algorithm (see also \cite{bonaventura17,edwards11}).


\begin{table}[!h]
\begin{center}
   \begin{tabular}{ c c | c c | c c }
\hline
\hline\noalign{\smallskip}
     $\dtnum{max}$ & 2000 & $\mathrm{R_{TOL}}$ & $10^{-6}$ & $\mathrm{A_{TOL}}$ & $10^{-3}$ \\ 
     $\mathrm{TOL}_N$ & $10^{-6}$ & $\kappa_N$ & 0.5 & $K$ & 10 \\ 
     $fac_{S1S2}$ & 0.3 & $fac_{\min}$ & 0.1 & $fac$ & $\sqrt[3]{0.25}$ \\
     $\mathrm{ratiomin}$ & 0.2 & $\mathrm{ratiomax}$ & 5.0 & $\mathrm{ratio}$ & $\sqrt[3]{0.25}$ \\
     $k_I$ & $\frac{1}{3}$ & $\epsilon$ & $10^{-10}$ & &\\
\noalign{\smallskip}\hline
\hline
   \end{tabular}
 \end{center}

\vspace{-3mm}
 \caption{Example 2. Adaptive timestep controller parameters.}\label{tab:dtAdaptParam1}
\end{table}

\begin{figure}[!h]
\begin{center}
\includegraphics[width=0.325\textwidth]{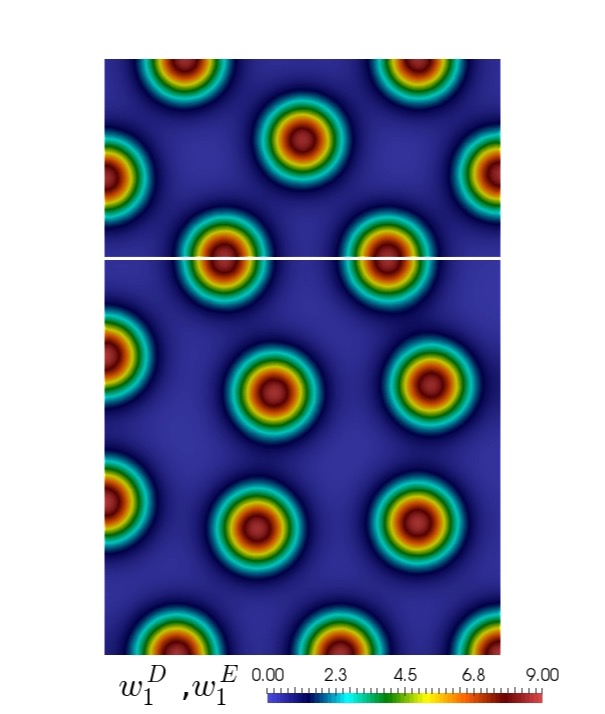}
\includegraphics[width=0.325\textwidth]{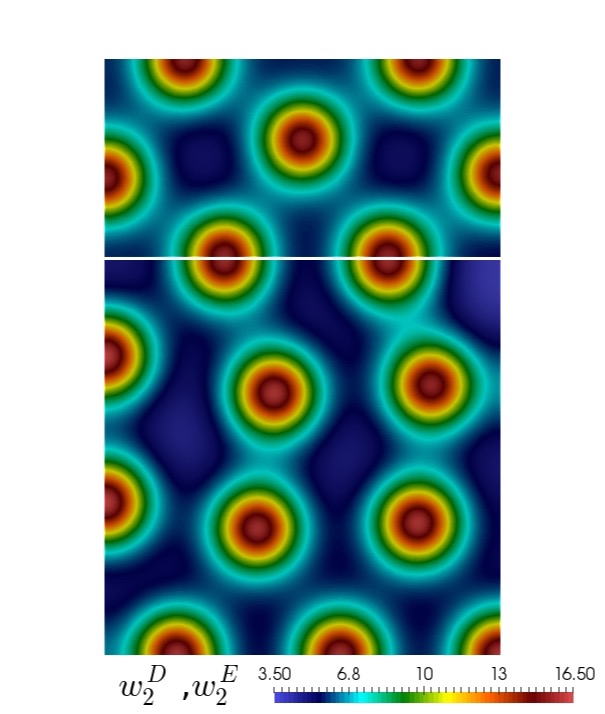}
\caption{Example 2. Species concentrations $w_1^*$, $w_2^*$ computed at $T = 2000$.}\label{fig:ex02-a} 
\end{center}
\end{figure}

\paragraph{Example 2: Cross-diffusion effects in a two-dimensional morphogen interface model.} In this example we focus our attention on the effect of cross-diffusion $\bM^*_{ij} \neq 0$ $\forall i,j$ in \eqref{eq:rd}. We consider again the Gierer-Meinhardt kinetics \eqref{reaction:GM}, now decoupled from the linear elasticity equations \eqref{eq:elasticity}. We impose the same reaction parameters on both sides of the interface $\rho_0 = 0$, $\rho_1 = \rho_2= \rho_4 = \rho_5= 1$, $\rho_3 = 0.35$, we set the diffusion coefficients $M_{11} = M_1 = 1$, $M_{22} = M_2 = 30$, and choose the acceleration transmission constants $K^\mathrm{D}  =  K^\mathrm{E} = 1$. The domains of interest are the 
blocks $\OmD=(0,50)^2$ and $\OmE = (0,50)\times(50,75)$, which we discretise into unstructured meshes with 16360 and 8032 triangular elements, respectively. Initial concentrations are set according to $$\tilde{w}_1 = \frac{\rho_2}{\rho_3}\biggl(\rho_0+\rho_1\frac{\rho_5}{\rho_4}\biggr) ,  \quad 
\tilde{w}_2 = \frac{\rho_4}{\rho_5} \tilde{w}_1^2 , \quad \bw^*_0(\bx) = \begin{pmatrix} \tilde{w}_1(1+\eta(\bx))\\ 
\tilde{w}_2\end{pmatrix},$$
where  $\eta(\bx)$ is a uniformly distributed random field with variance $10^{-3}$. We employ the adaptive time-stepping method described in Algorithm \ref{algo:TRBDF2} with controller parameters defined as in Table \ref{tab:dtAdaptParam1}. We simulate the process until $T=2000$, and display the resulting numerical solutions on Figure \ref{fig:ex02-a}.

\begin{figure}[!h]
\begin{center}
\includegraphics[width=0.85\textwidth]{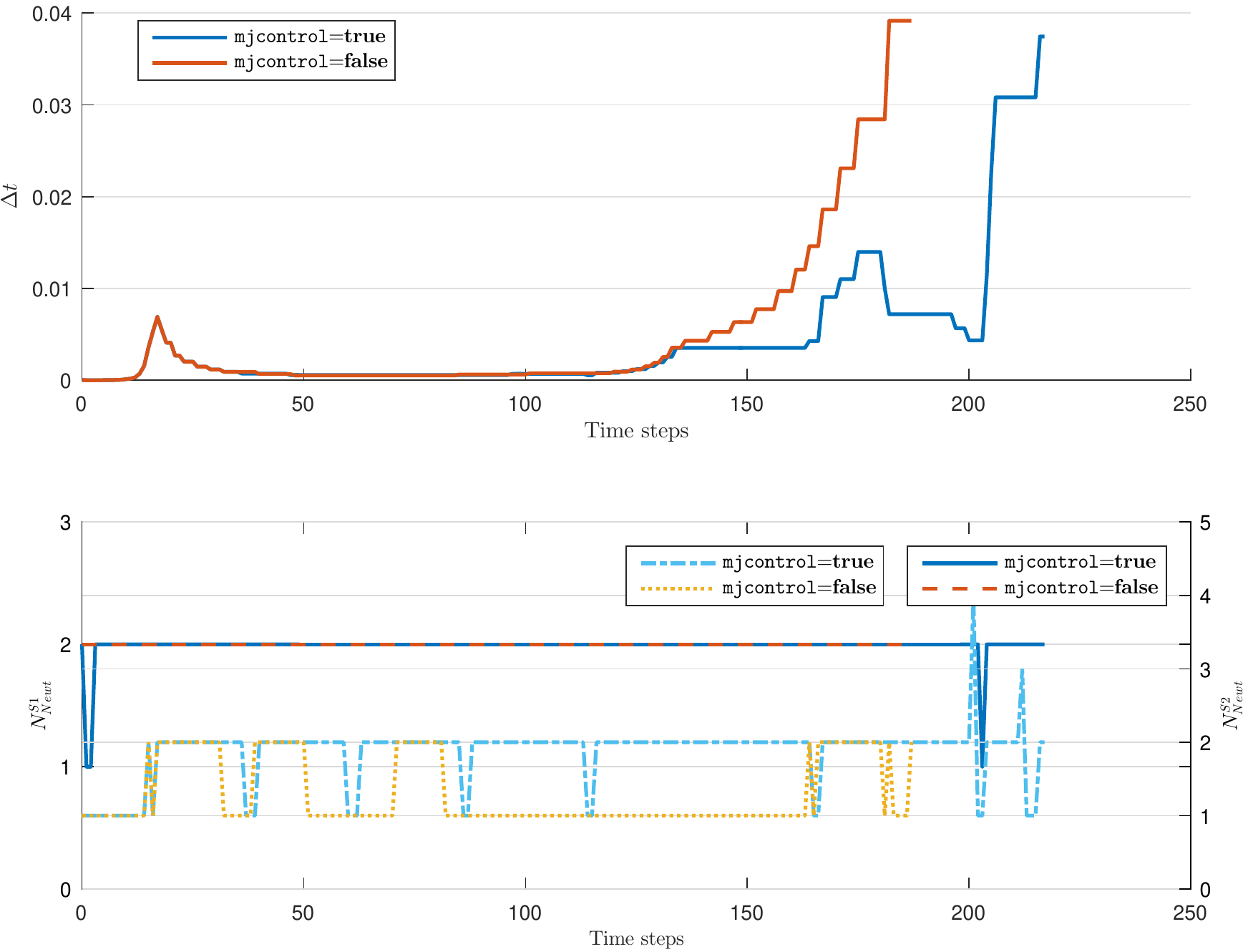}
\caption{Example 2. Time step evolution and number of Newton iterations needed 
to converge on stages \textsc{S1} and \textsc{S2}, and setting \texttt{mjcontrol} on and off.}\label{fig:ex02-b} 
\end{center}
\end{figure}

In Figure~\ref{fig:ex02-b} we observe that the adaptive time stepping behaves differently depending on whether one reuses or not the Newton iteration matrices over the implicit stages (this feature is turned on/off by the parameter \texttt{mjcontrol}). For each case, we observe an increase of $\dtnum{}{}$ as the system reaches a 
steady state. The spatial patterns generated by the two approaches practically coincide and 
therefore are not shown. However, we observe that the case $\texttt{mjcontrol}=\True$ does not show 
a marked increase of the timestep when the system reaches stationary patterns. Additionally, we 
observe that reusing the iteration matrix leads to a slight increase in the number of Newton 
iterations in the \textsc{S2} stage.

\begin{figure}[t]
\begin{center}
\includegraphics[width=0.325\textwidth]{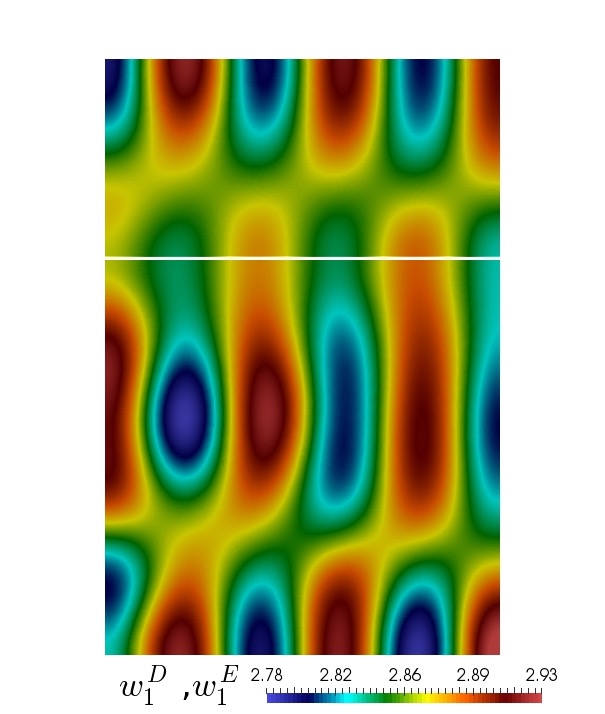}
\includegraphics[width=0.325\textwidth]{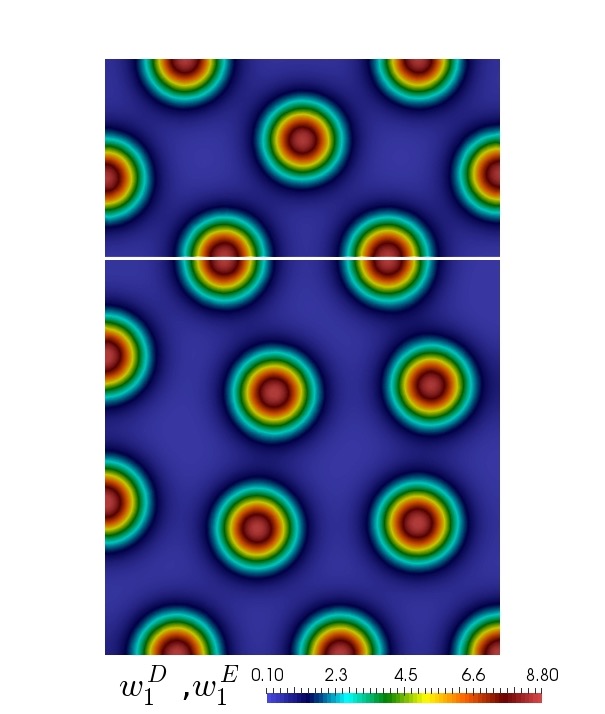}
\includegraphics[width=0.325\textwidth]{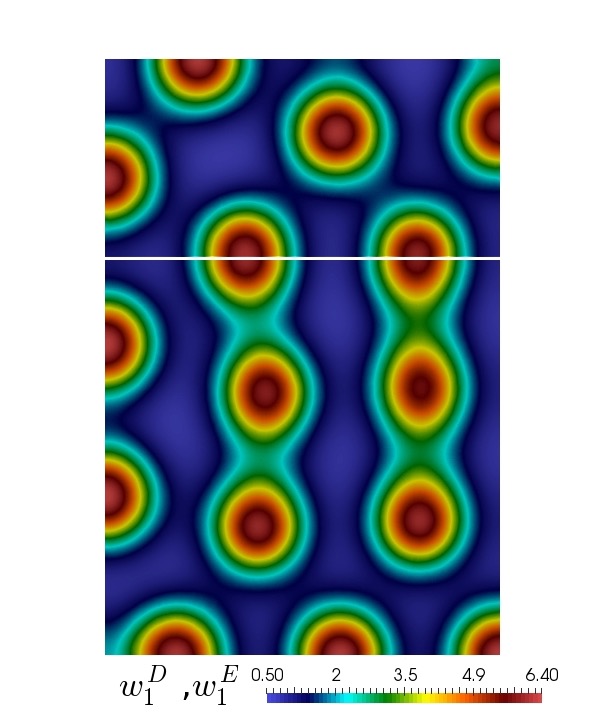}
\includegraphics[width=0.325\textwidth]{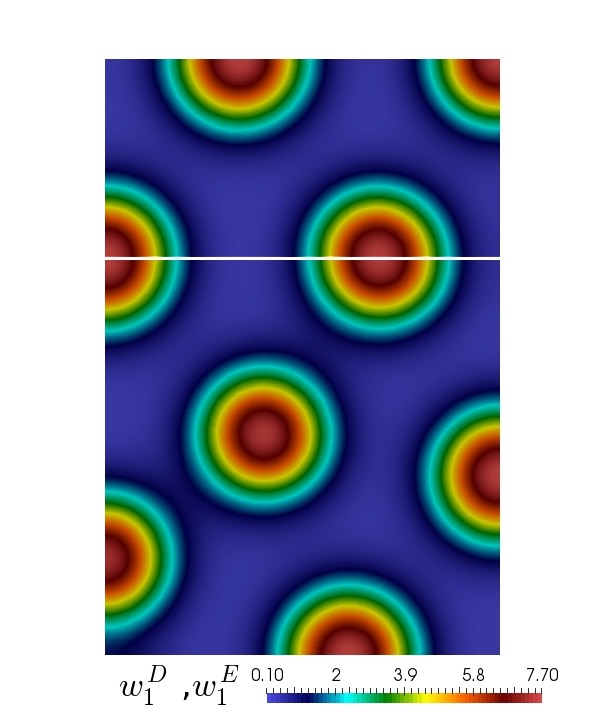}
\includegraphics[width=0.325\textwidth]{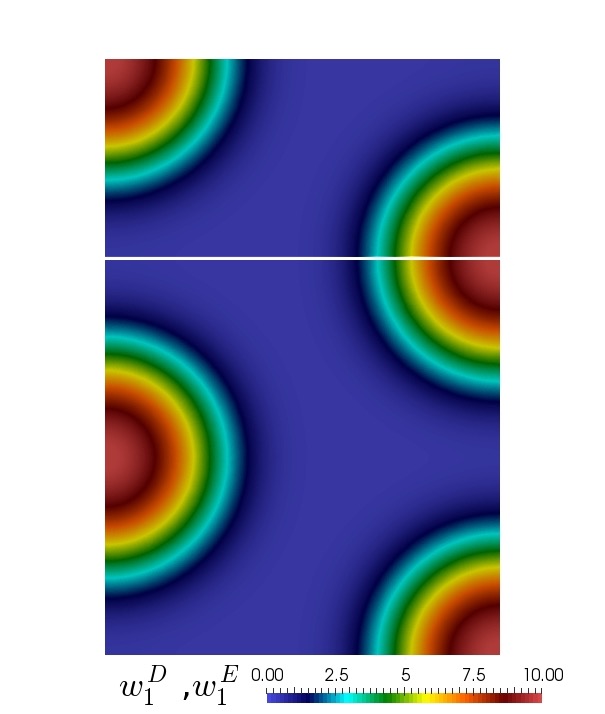}
\includegraphics[width=0.325\textwidth]{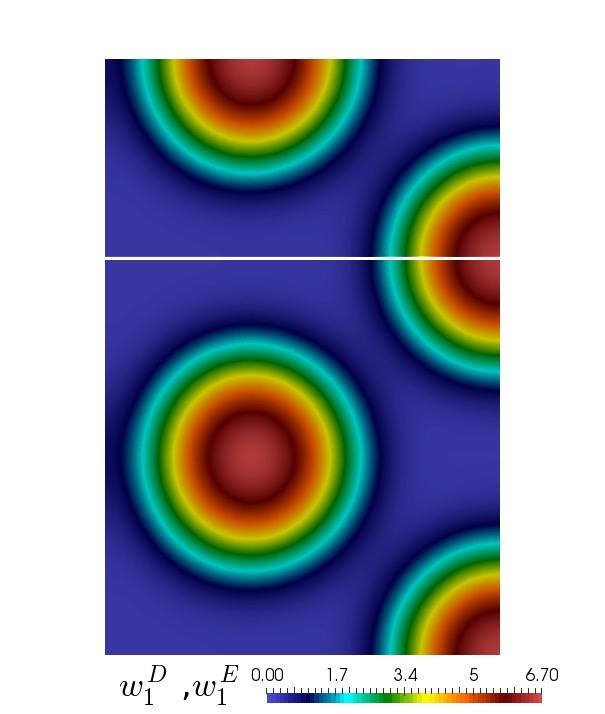}
\caption{Example 2, 2D case. 
Concentration of species $w_1^*$. Linear cross-diffusion (top row) using 
$(M_{12}^*,M_{21}^*)=(0.5,0)$, $(0,0.5)$ and $(0.35,0.05)$ (left, centre, right, respectively). 
The bottom panels show the case of nonlinear cross-diffusion using 
$(\eta_{11}^{*,w_1},\eta_{11}^{*,w_2},\eta_{22}^{*,w_1},\eta_{22}^{*,w_2}=(10^{-1},10^{-1},10^{-1},10^{-1})$, $(1,0,0,1)$ and $(1,0,1,0)$ (left, centre, right, respectively).}
\label{fig:ex02-c} 
\end{center}
\end{figure}

\begin{figure}[!h]
\begin{center}
\includegraphics[width=0.415\textwidth, clip=true, trim={0cm 0cm 0cm 0cm}]{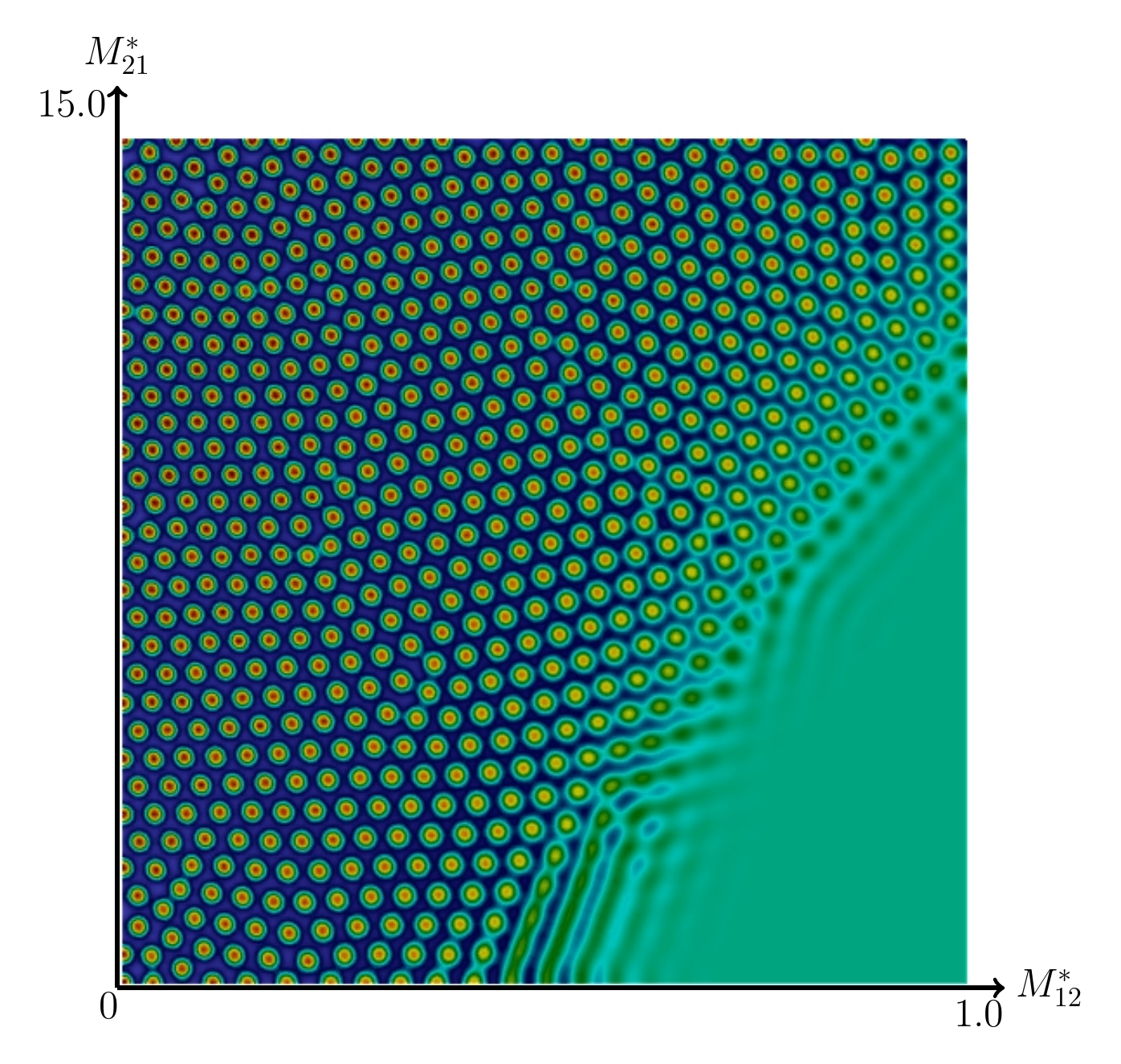}
\includegraphics[width=0.4\textwidth, clip=true, trim={0cm 0cm 0cm 0cm}]{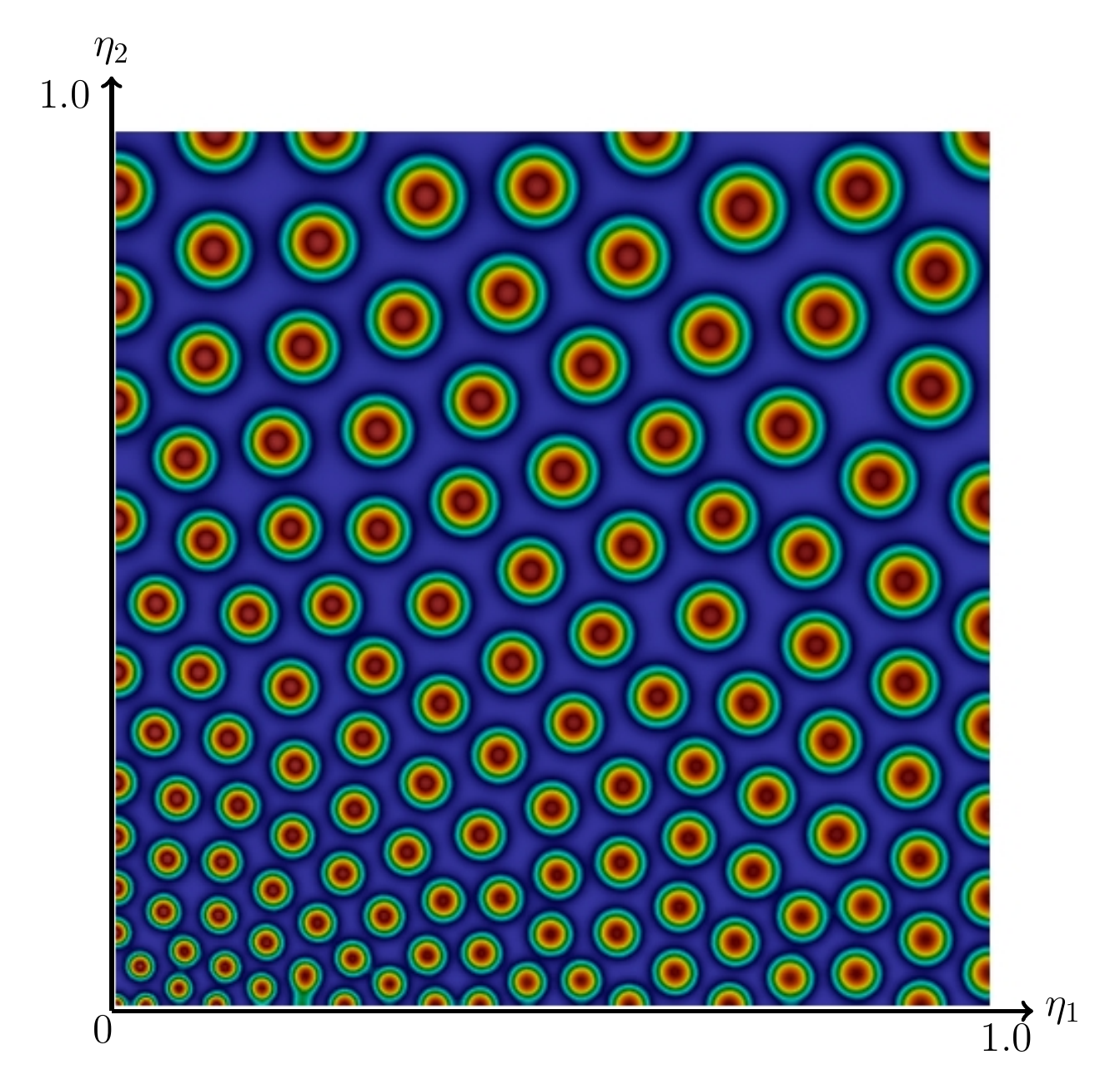}
\caption{Example 2. Space parameter plots produced with linear (left) and nonlinear (right) 
cross-diffusion terms.}\label{fig:ex02-d} 
\end{center}
\end{figure}

\begin{figure}[t]
\begin{center}
\includegraphics[width=0.325\textwidth]{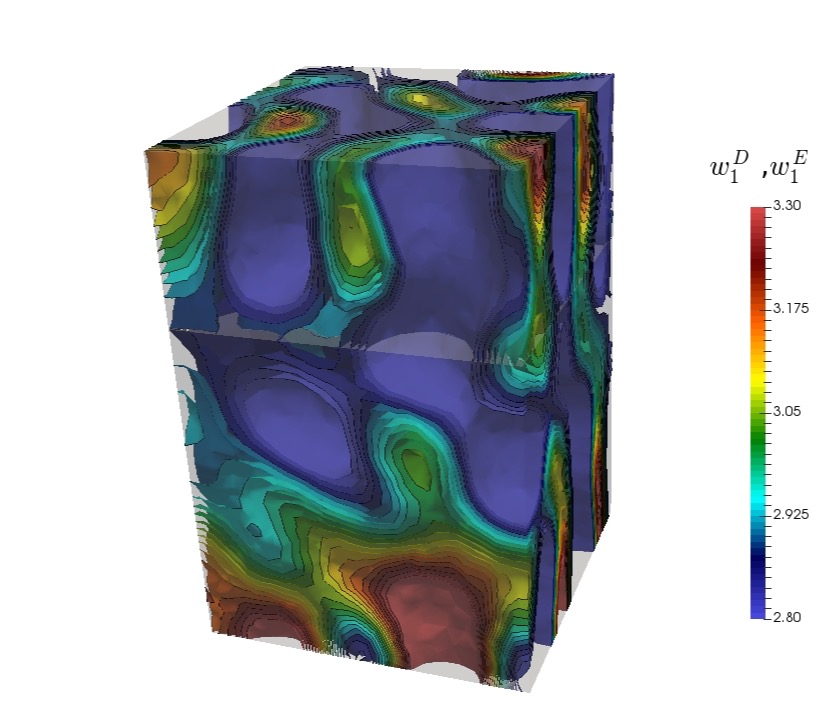}
\includegraphics[width=0.325\textwidth]{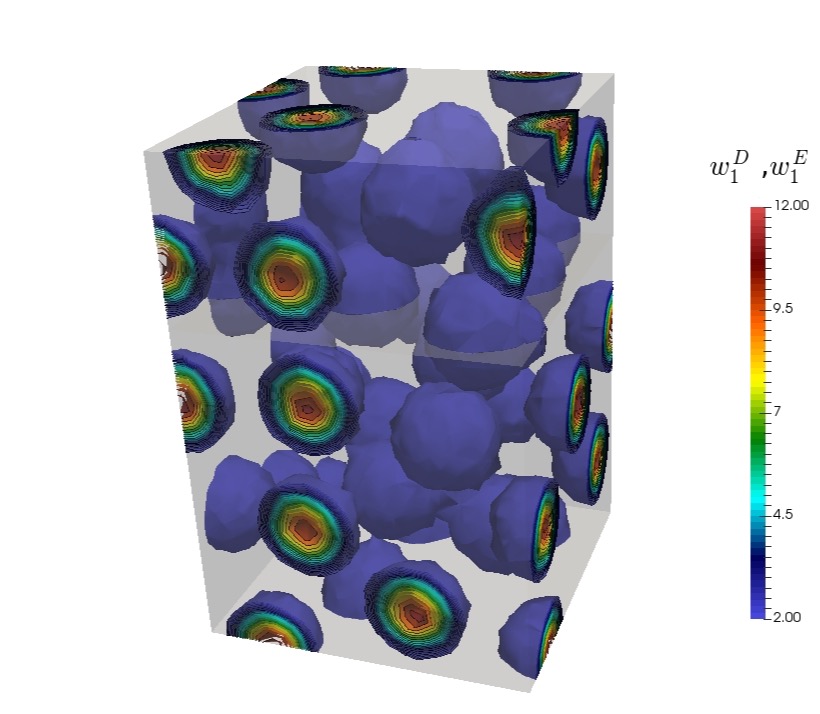}
\includegraphics[width=0.325\textwidth]{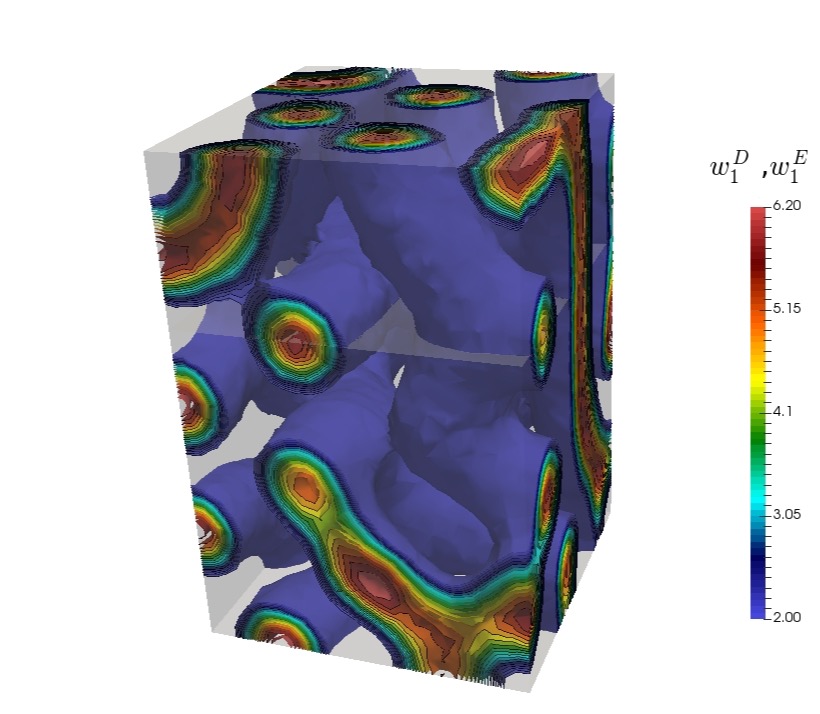}
\includegraphics[width=0.325\textwidth]{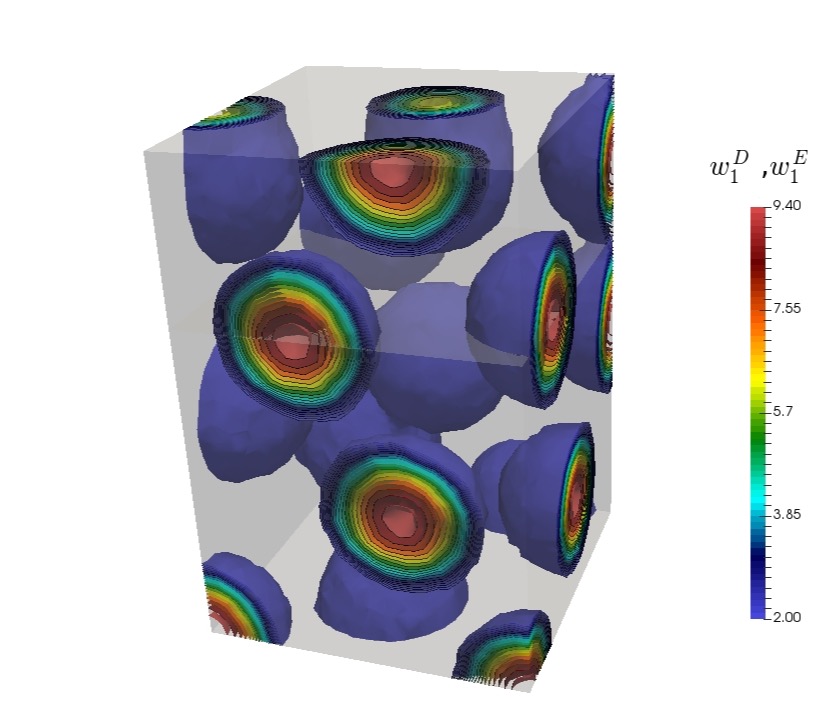}
\includegraphics[width=0.325\textwidth]{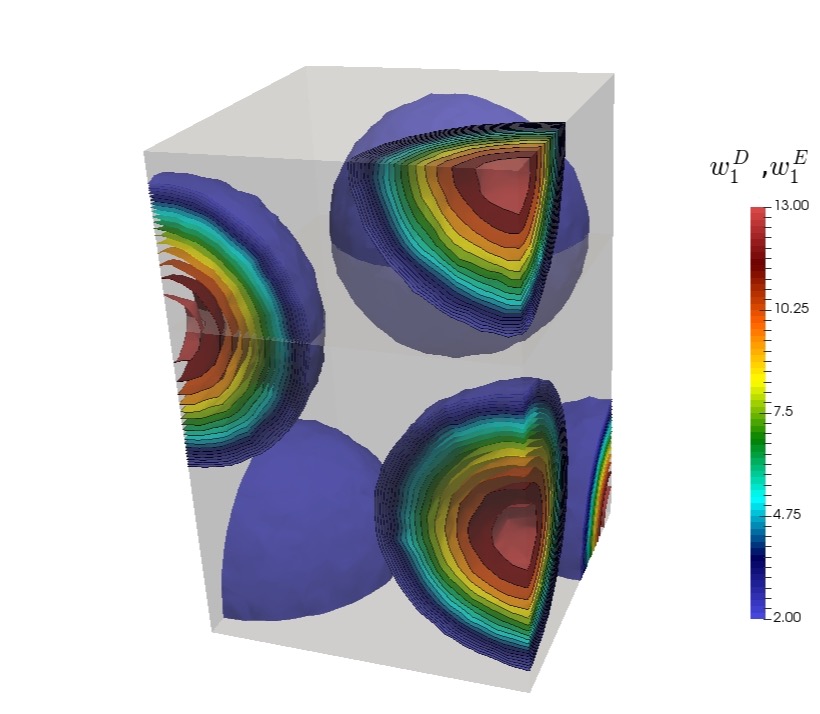}
\includegraphics[width=0.325\textwidth]{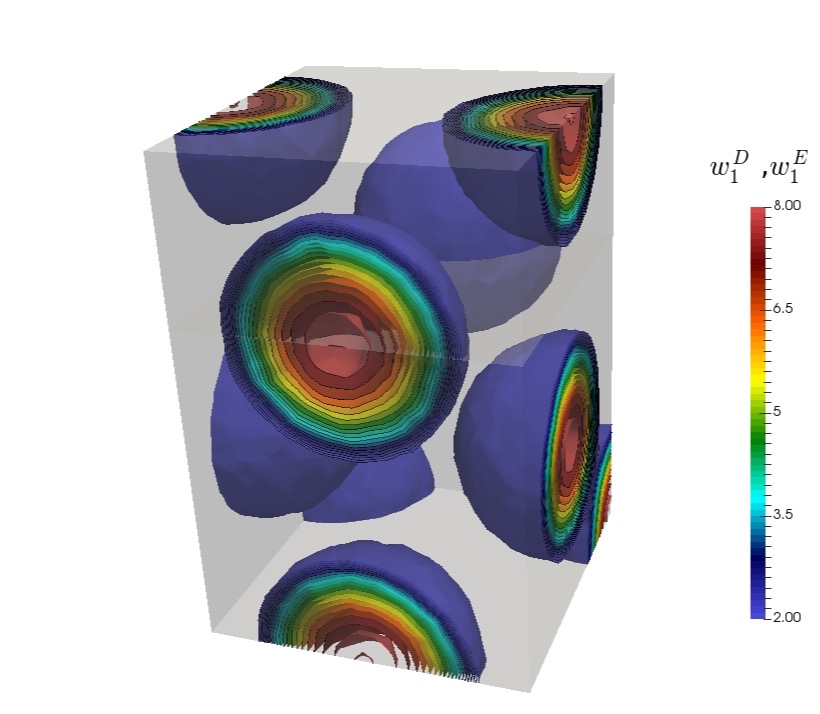}
\caption{Example 2, 3D case. 
Concentration of species $w_1^*$. Linear cross-diffusion (top row) using 
$(M_{12}^*,M_{21}^*)=(0.5,0)$, $(0,0.5)$ and $(0.35,0.05)$ (left, centre, right, respectively). 
The bottom panels show the case of nonlinear cross-diffusion using 
$(\eta_{11}^{*,w_1},\eta_{11}^{*,w_2},\eta_{22}^{*,w_1},\eta_{22}^{*,w_2})=(10^{-1},10^{-1},10^{-1},10^{-1})$, $(1,0,0,1)$ and $(1,0,1,0)$ (left, centre, right, respectively).}
\label{fig:ex02-e} 
\end{center}
\end{figure}

In addition, we test our algorithm for both linear and nonlinear cross-diffusion 
terms defined as 
$$\bM_{lin}^*(\bw^*) = \begin{pmatrix} 
M^*_{11} & M^*_{12} \\
M^*_{21} & M^*_{22} \end{pmatrix}
$$
or
$$\bM_{\mathrm{nl}}^*(\bw^*) = \begin{pmatrix} 
M^*_{11}(1 + \eta_{11}^{*,w_1}w_1^* + \eta_{11}^{*,w_2}w_2^*) & 0  \\
0 & M^*_{22}(1 + \eta_{22}^{*,w_1}w_1^* + \eta_{22}^{*,w_2}w_2^*) \end{pmatrix},
$$
where $M_{ij}^*$, $\eta_{ij}^{*,w_k}$ ($i,j,k=\{1,2\}$) are real constants. For these 
simulations we impose $M^*_{11}=M^*_1$, $M^*_{22}=M^*_2$. As illustrated in Figure \ref{fig:ex02-c}, the final patterns adopt different structures according to the cross-diffusion matrix employed. While the specific nonlinear cross-diffusion matrix modifies the size of $w_i$-concentration spots without changing its structure, the linear one (and depending on the parameter values) can change from spots to striped patterns. Furthermore, nonlinear cross-diffusion induces a more visible impact on the distance between high concentration regions, when compared to the linear case. Space-parameter plots illustrate this phenomenon in Figure \ref{fig:ex02-d}. For both cross-diffusion models, we vary spatially the values $[M_{12}^*$, $M_{21}^*] \in [0,1]\times[0,15]$, and $\eta_{11}^{*,w_1}=\eta_{22}^{*,w_1}=\eta_1$ and $\eta_{11}^{*,w_2}=\eta_{22}^{*,w_2}=\eta_2$ from 0 to 1; always checking that the cross-diffusion matrices remain positive definite. 
The linear cross-diffusion particularly affects the type of patterns that are generated but does not significantly change length scale of the pattern. Nonlinear diffusion, on the other hand, affects directly the pattern length scale without altering the general pattern structure produced by the system. These observations are confirmed in 3D domains as well (cf. \ref{fig:ex02-e} ).

\begin{figure}[t]
\begin{center}
\includegraphics[width=0.325\textwidth]{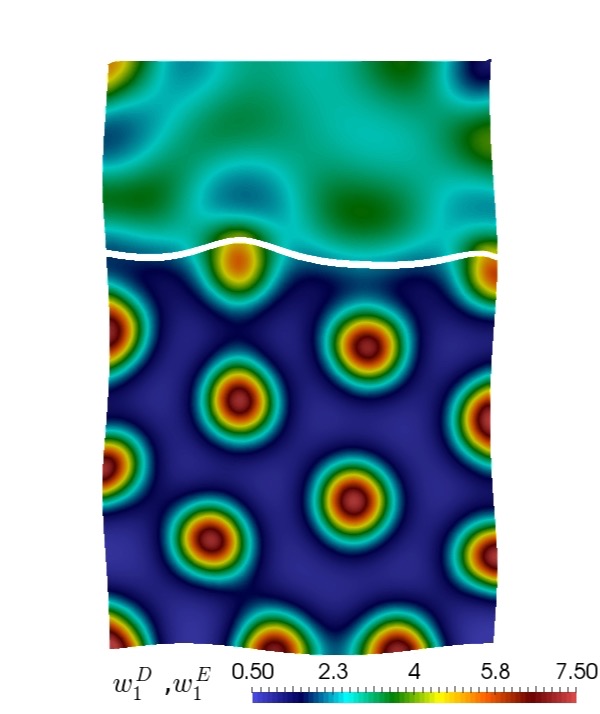}
\includegraphics[width=0.325\textwidth]{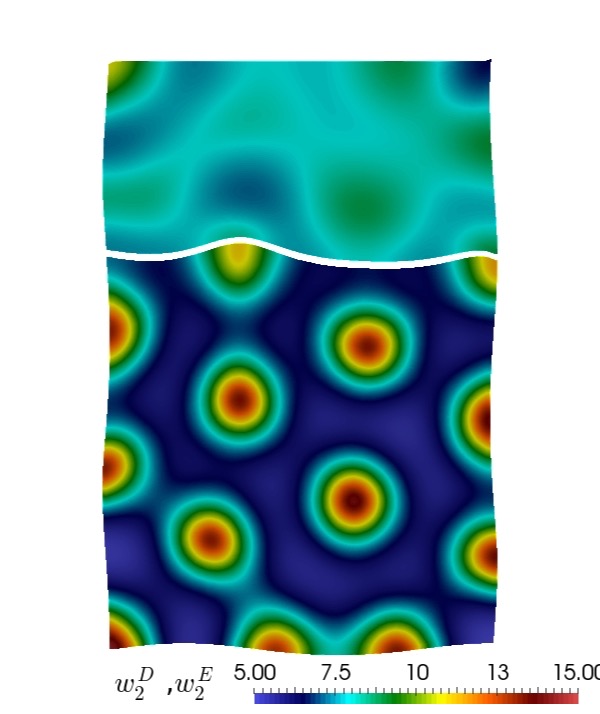}
\includegraphics[width=0.335\textwidth]{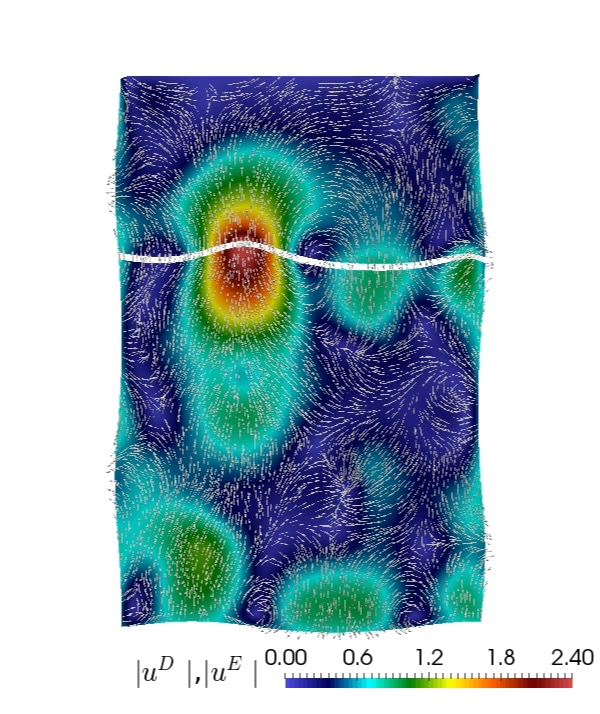}
\caption{Example 3. Morphogen concentrations $w_1^*$, $w_2^*$ (left and centre)
and displacement magnitude (right), computed at $T = 2000$ and plotted 
on the deformed configuration. }\label{fig:ex03-a} 
\end{center}
\end{figure}

\paragraph{Example 3: A simplified activator-inhibitor model on an elastic 2D domain.} 
With the objective of testing the time-adaptive algorithms  
now applied to the mechanochemical coupled problem, we add the elasticity components to the reaction kinetics \eqref{reaction:GM} of Example 2 (with the same parameters and mesh). The Young moduli and Poisson ratios defining the material properties of the medium are now $E^\mathrm{D} = 1000$, $E^\mathrm{E} = 250$, $\nu^\mathrm{D} = 0.475$, $\nu^\mathrm{E} = 0.3$. Displacements 
are prescribed on the boundary $\partial\Omega\setminus\Gamma$ according to 
 $\bu = (\frac{1}{2}\cos(5\pi(x-y)/2),\frac{3}{4}\sin(5\pi (x+y)/2) )^T$. Additionally, we impose $c_f^{\mathrm{D}}=150$, $c_f^{\mathrm{E}}=20$ and $c_g^*=1$. The medium is initially at rest and stress-free (zero displacement and pressure on the whole domain) and we re-use the initial conditions of the morphogen concentrations
defined in the previous test. We use Algorithm \ref{algo:TRBDF2} with controller parameters defined again as in Table \ref{tab:dtAdaptParam1}, except for $\mathrm{A_{TOL}}=10^{-6}$ and $\mathrm{TOL}_N = \mathrm{R_{TOL}}=10^{-4}$. The process is simulated until $T=2000$. 

Figure \ref{fig:ex03-a} presents the final patterns generated by the morphogen 
concentrations and also the zones of high displacement. On the bottom layer 
the pattern is similar to the one in Figure \ref{fig:ex02-a} (in terms of magnitude and 
spatial structure). For the layer with smaller Poisson ratio, the patterns tend to homogenise 
and eventually disappear. This phenomenon highlights the effects of mechanical properties 
on the morphogen dynamics, and  we also note that deformations are larger near the interface 
and close to the bottom boundary of $\OmD$. 

\begin{figure}[t!]
\begin{center}
\includegraphics[width=0.24\textwidth]{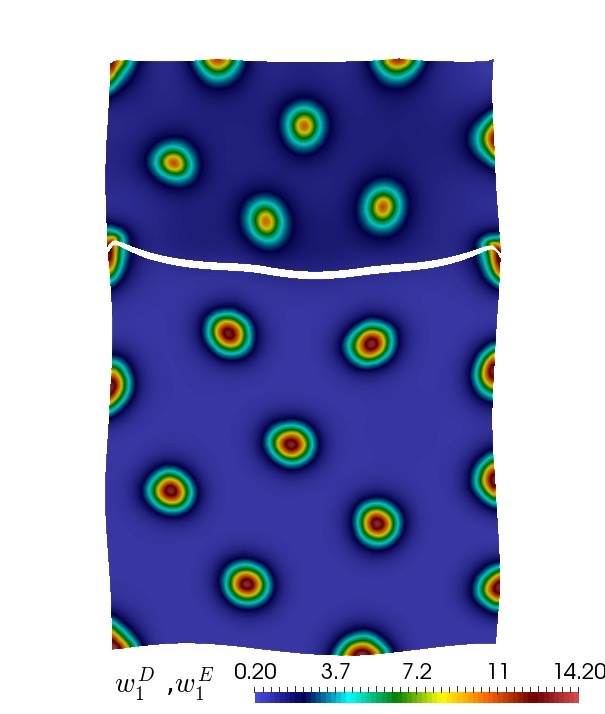}
\includegraphics[width=0.24\textwidth]{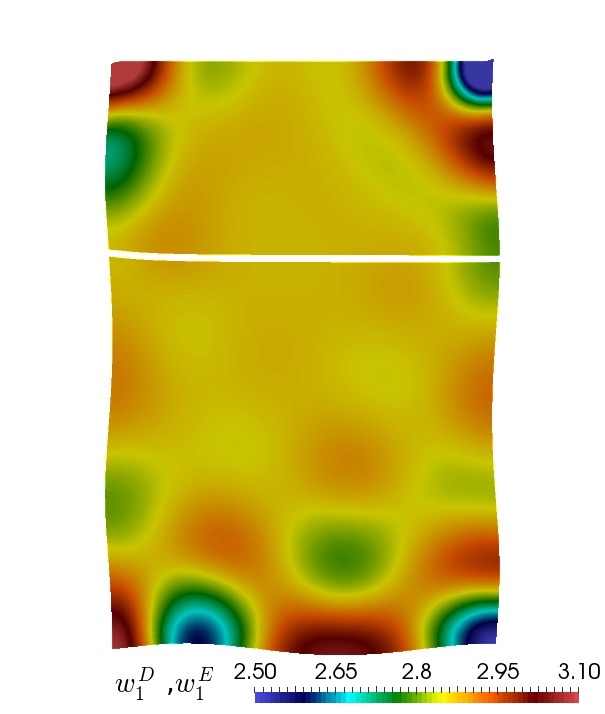}
\includegraphics[width=0.24\textwidth]{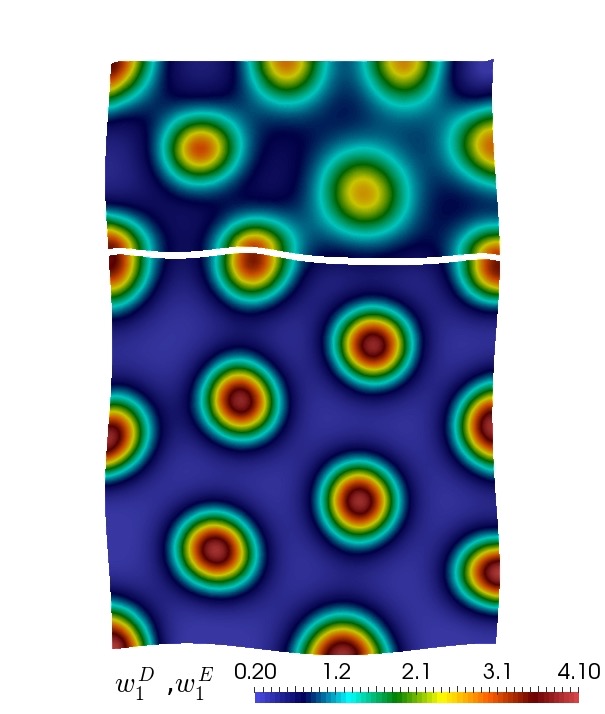}
\includegraphics[width=0.24\textwidth]{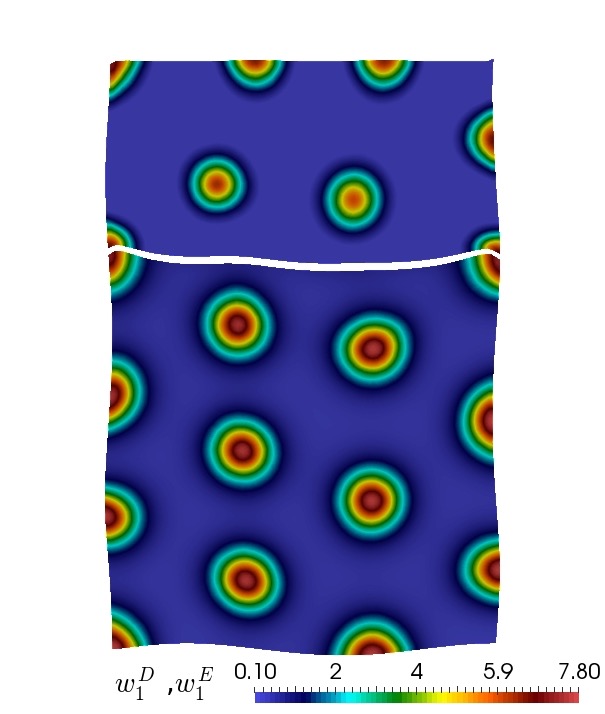}\\
\includegraphics[width=0.24\textwidth]{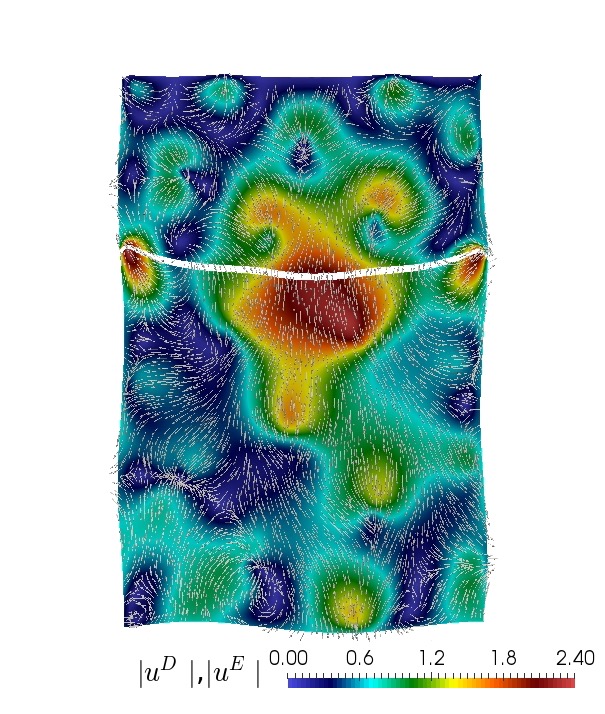}
\includegraphics[width=0.24\textwidth]{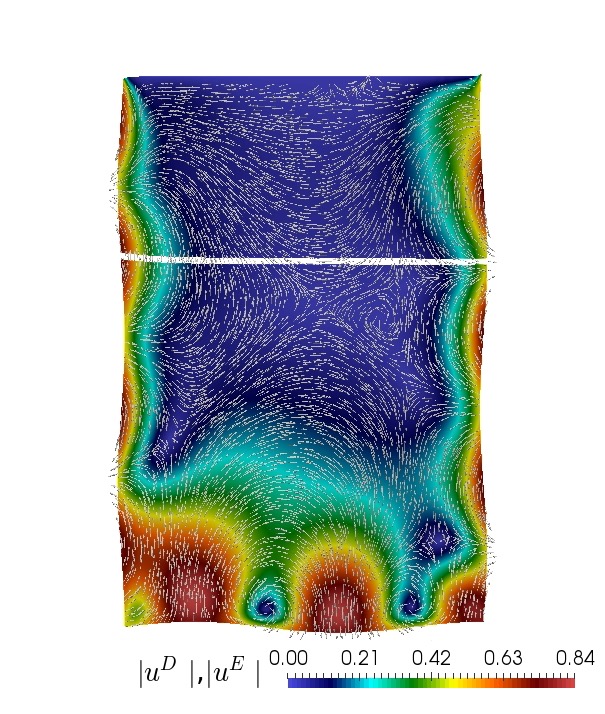}
\includegraphics[width=0.24\textwidth]{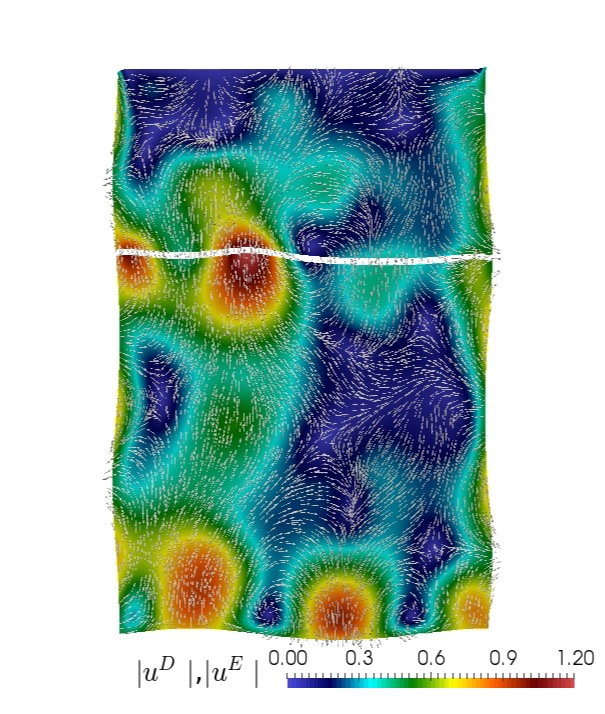}
\includegraphics[width=0.24\textwidth]{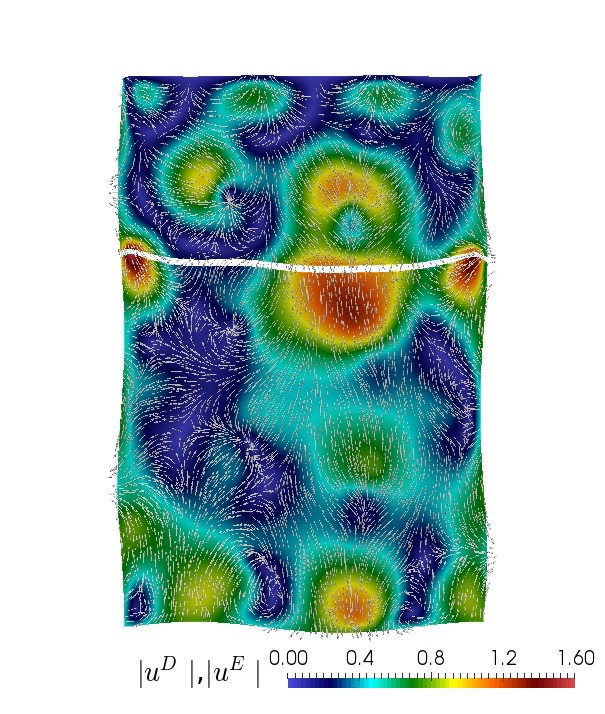}\\
\includegraphics[width=0.24\textwidth]{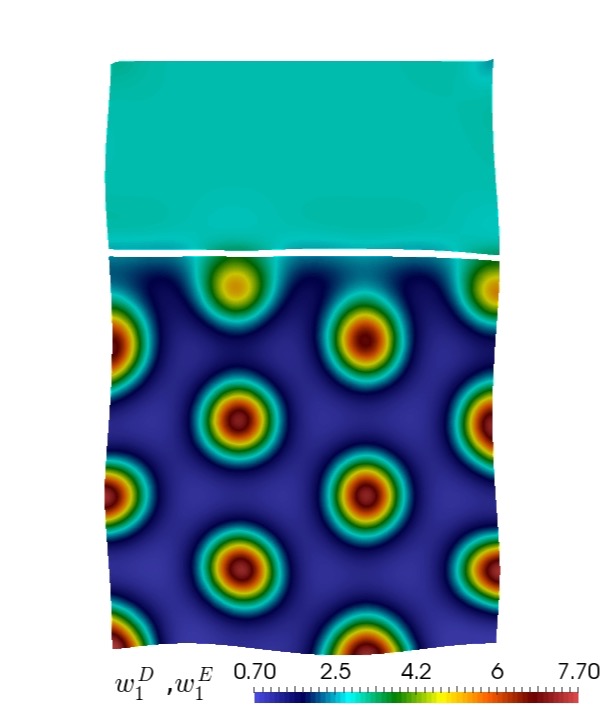}
\includegraphics[width=0.24\textwidth]{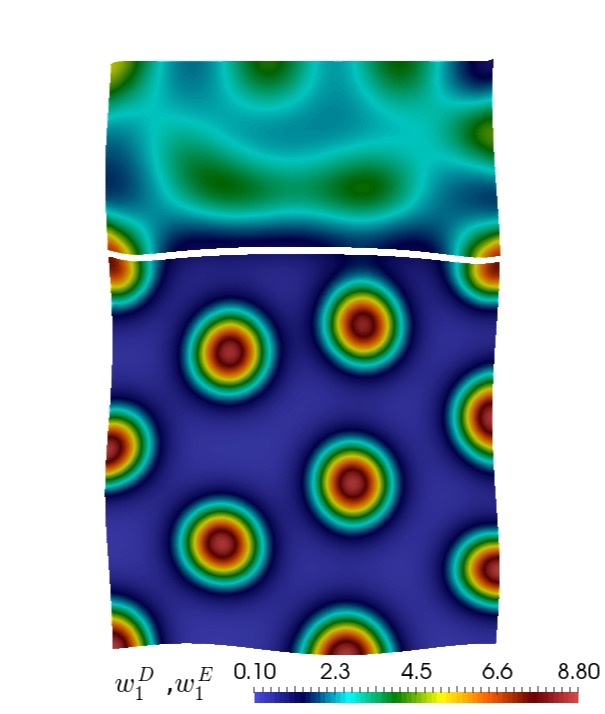}
\includegraphics[width=0.24\textwidth]{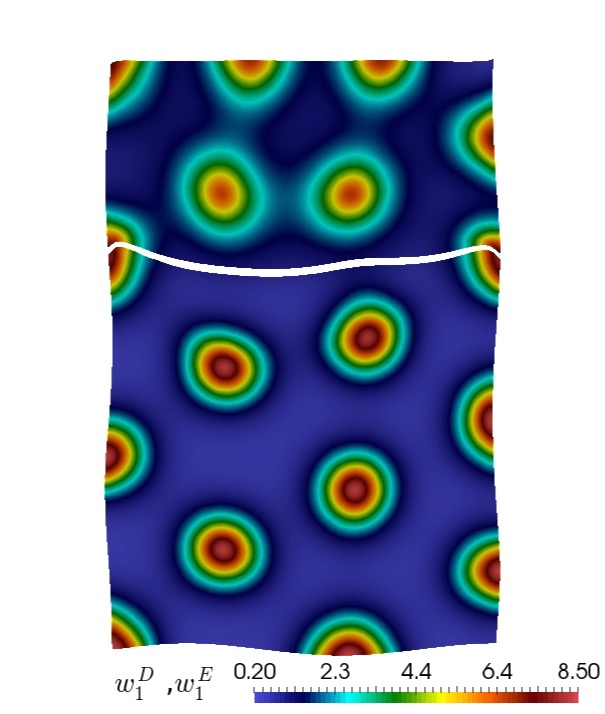}
\includegraphics[width=0.24\textwidth]{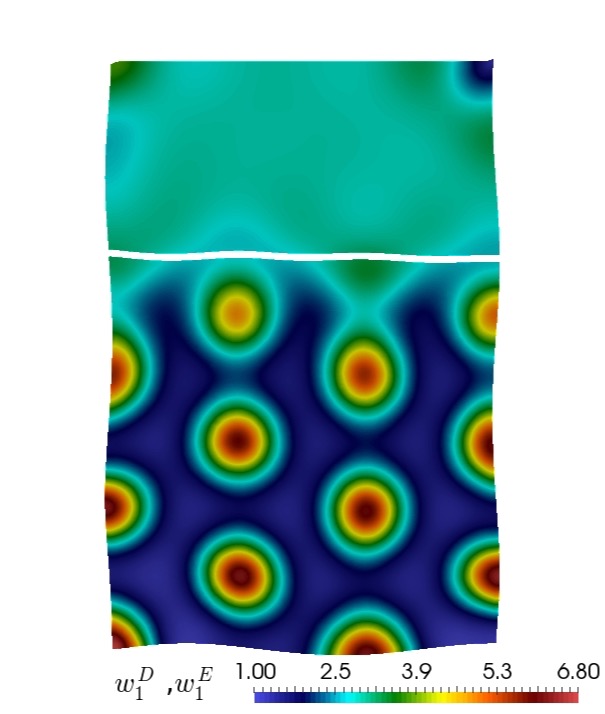}\\
\includegraphics[width=0.24\textwidth]{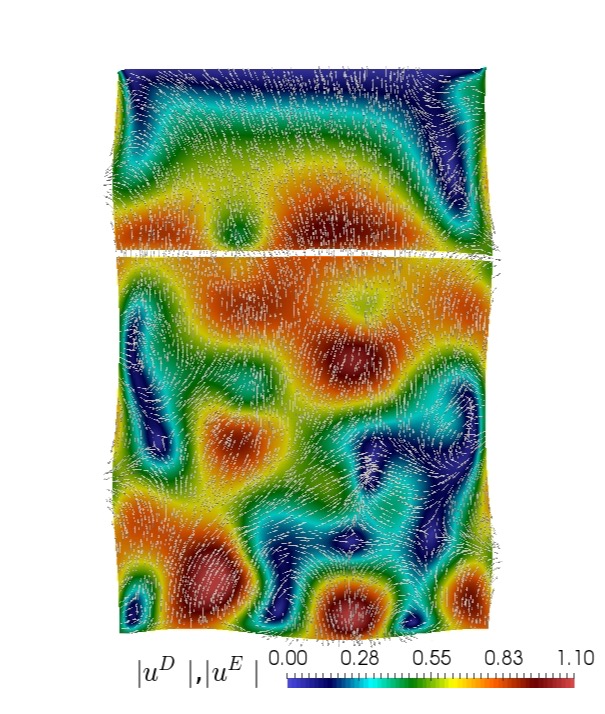}
\includegraphics[width=0.24\textwidth]{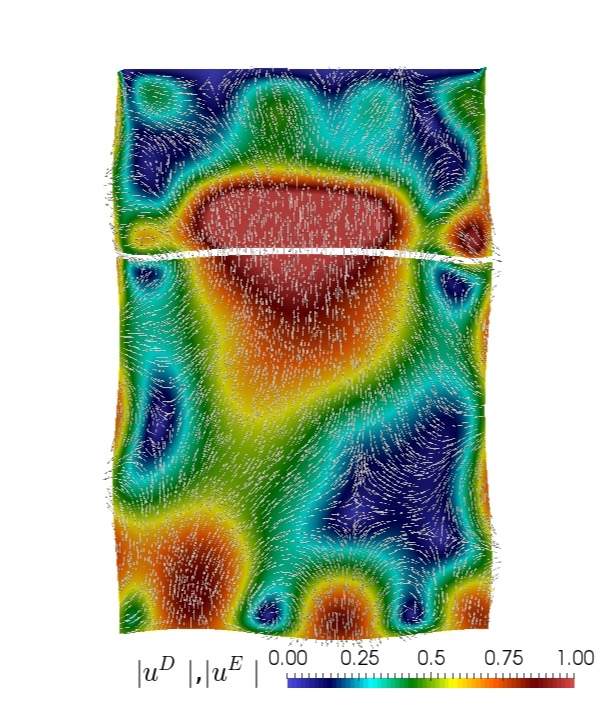}
\includegraphics[width=0.24\textwidth]{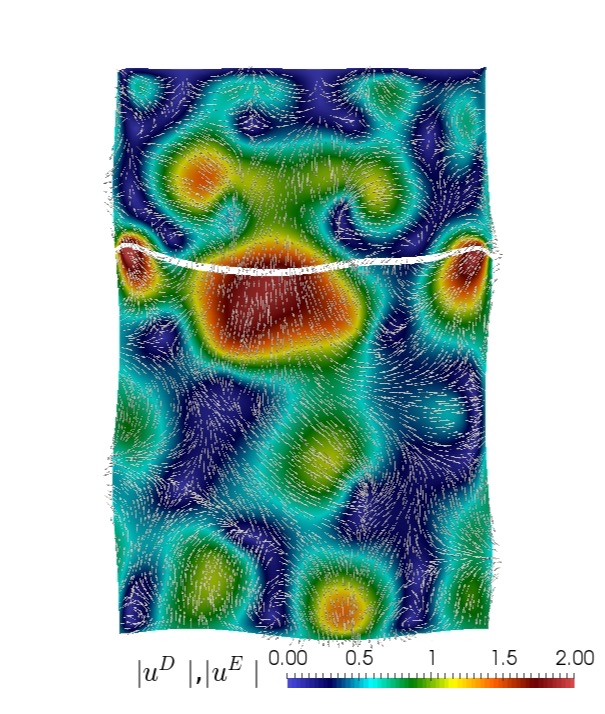}
\includegraphics[width=0.24\textwidth]{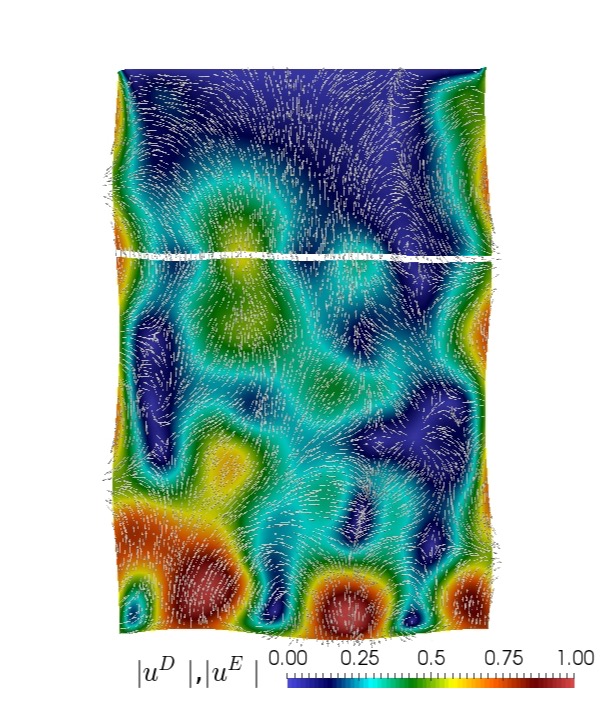}
\caption{Example 3. Concentration of species $w_1^*$ and displacement magnitude $|u^*|$. Influence of ADR parameters ($1^{\text{st}}$ \& $2^{\text{nd}}$ rows) using $M_1^*=0.5$, $M_1^*=1.5$, $\rho_2^*=0.5$ and $\rho_3^*=0.45$ (from left to right respectively); and of coupling  parameters ($3^{\text{rd}}$ \& $4^{\text{th}}$ rows) using $c_f^*=150$, $c_f^*=20$, $c_g^*=0.5$, $c_g^*=1.5$ (from left to right respectively). }\label{fig:ex03-b} 
\end{center}
\end{figure}

We next test the influence of the reaction and mechanical coupling parameters $(c_f^*,c_g^*)$ on the final pattern, and typical results are depicted in Figure \ref{fig:ex03-b}. Perturbations of the ADR 
coefficients are responsible for considerable variations on the species concentration spatial distribution as well as on the deformation profiles produced by the model. For instance, we observe that the size and shape of morphogens concentration patterns induce differences on the deformation of the interface $\Sigma$. Conversely, differences on the mechanical response of each subdomain also directly impact on the chemical species distribution, at the interface. Surprisingly, varying the coupling parameters does not substantially impact the length scale but rather make these patterns to disappear altogether. 

\begin{figure}[!h]
\begin{center}
\includegraphics[width=0.475\textwidth]{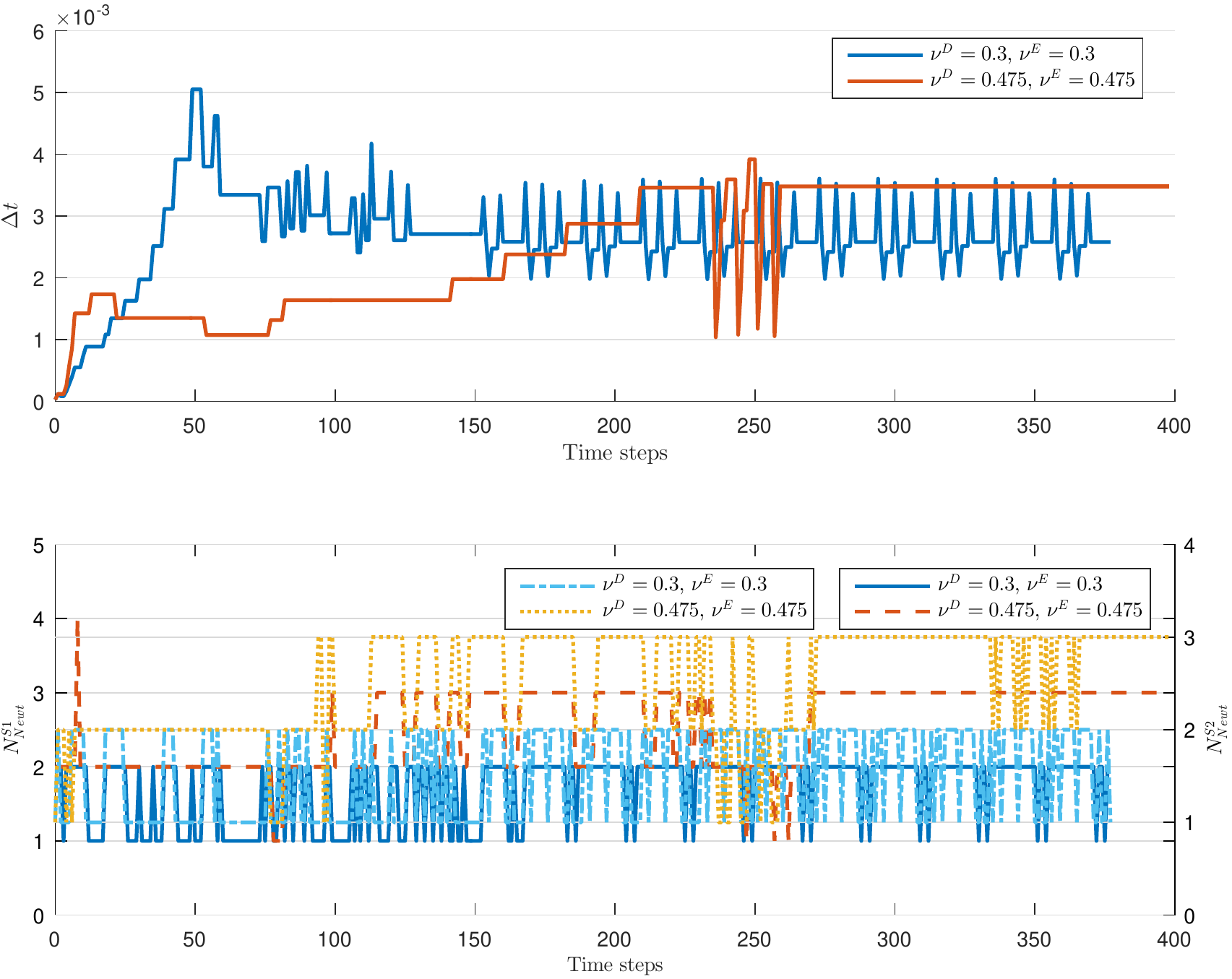}
\includegraphics[width=0.475\textwidth]{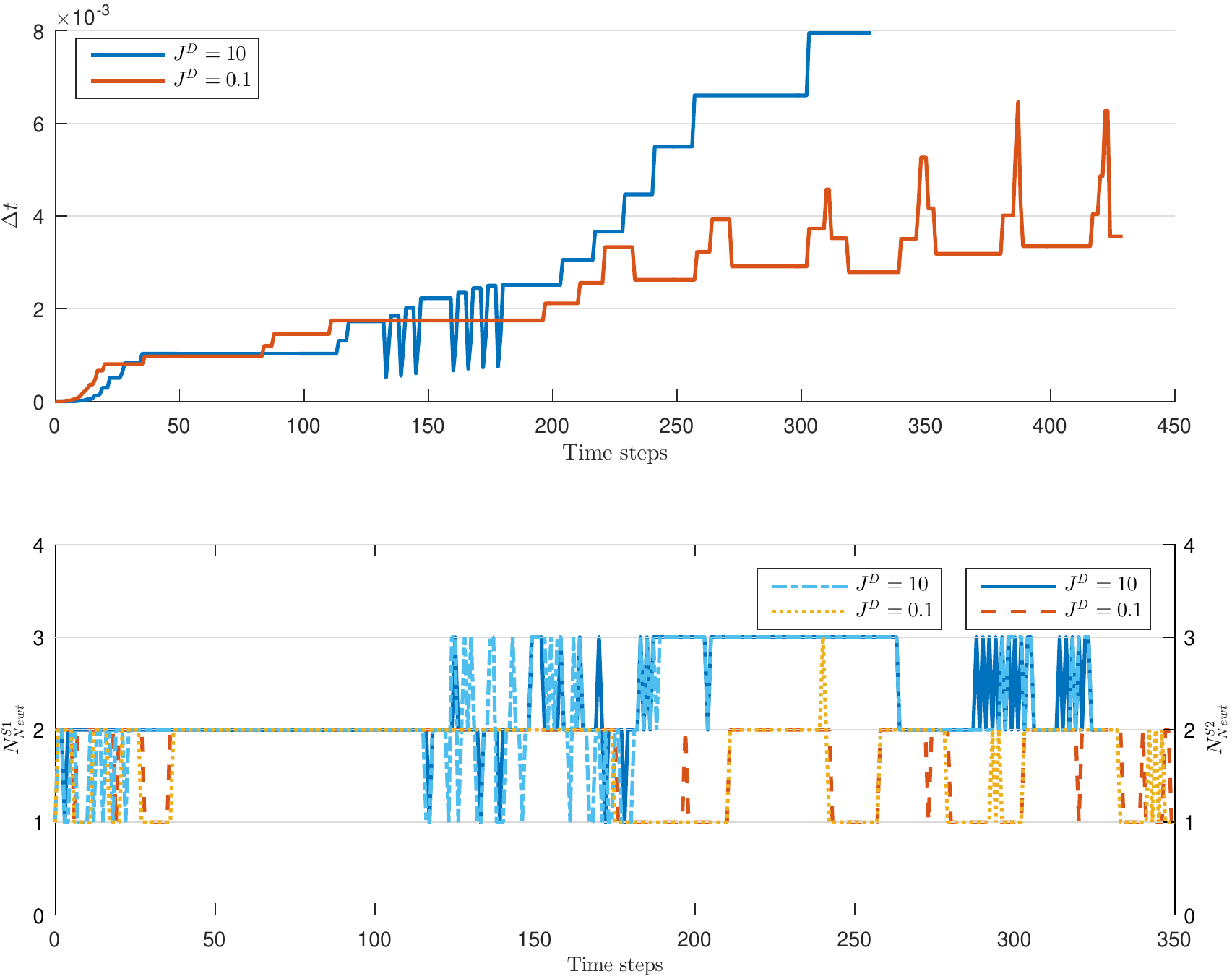}
\includegraphics[width=0.245\textwidth]{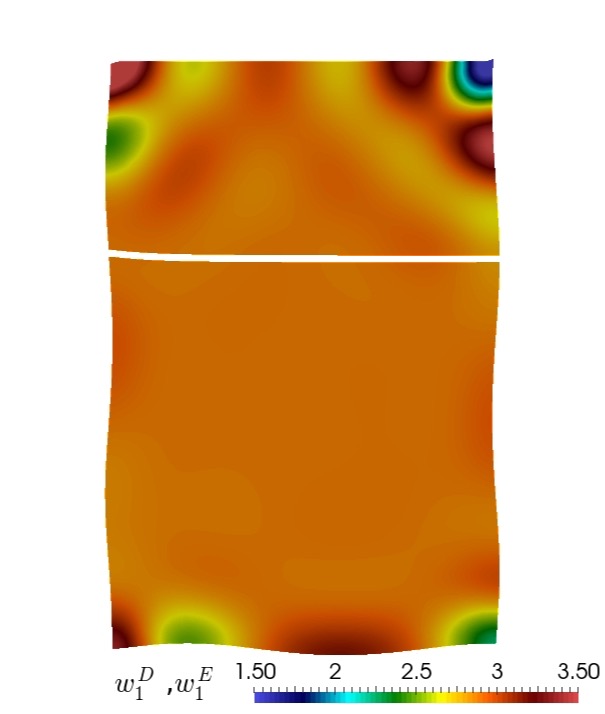}
\includegraphics[width=0.245\textwidth]{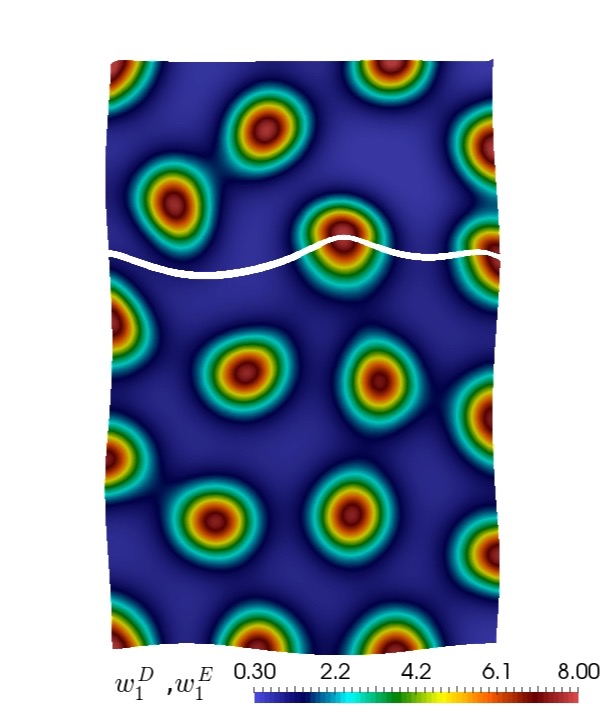}
\includegraphics[width=0.245\textwidth]{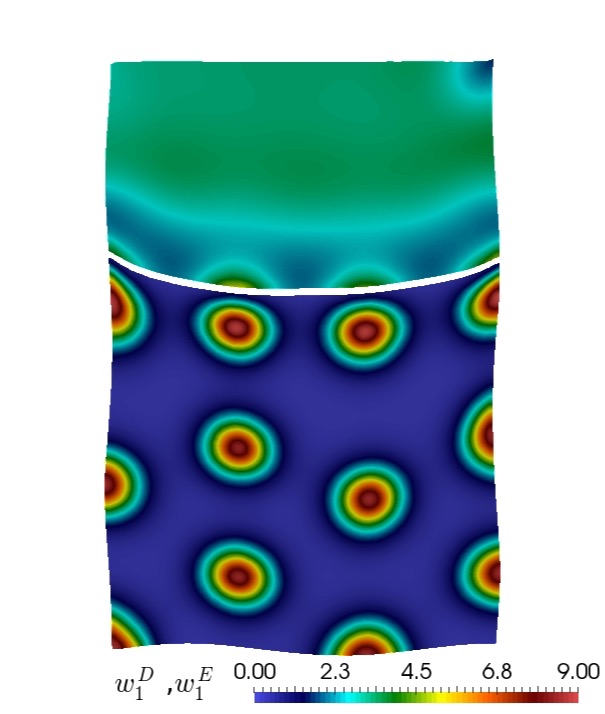}
\includegraphics[width=0.245\textwidth]{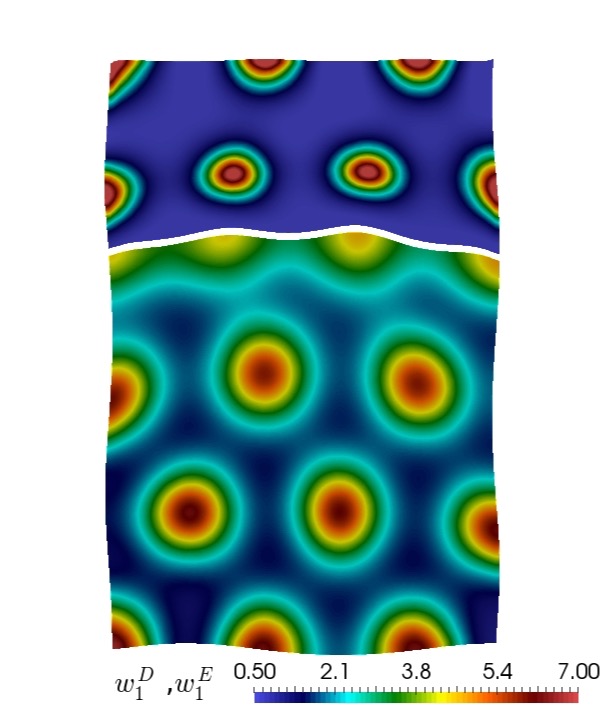}
\caption{Example 3. Top: Timestep evolution and Newton iteration count at stages \textsc{S1} and \textsc{S2} for different mechanical properties and acceleration constants. Species concentration $w_1^*$ using $\nu^*=0.3$, $\nu^*=0.475$, $J^D=10$ and $J^D = 0.1$ (bottom, from left to right, respectively).}\label{fig:ex03-c} 
\end{center}
\end{figure}

We also inspected the effect of the acceleration constants of the ADR interface conditions 
together with the degree of material heterogeneity (how different the Lam\'e parameters are 
in each subdomain) on the timestep evolution and on the Newton iteration count. 
For these simulations we fix \texttt{mjcontrol} = \True in order to inspect also the efficiency 
on the update of the Jacobian matrix throughout the simulation. Some of these results are collected in 
Figure \ref{fig:ex03-c}. The value $N^{SX}_{Newt}$ represents the number of iterations corresponding to 
\emph{accepted} timesteps. For this case we observe that $\dtnum{}{}$ is of the order of $10^{-3}$ and 
progressively increases with the simulation time and the consolidation of the spatial patterns. 
The algorithm readily detects when the process approaches a steady state, leading to a fast increase of the timestep. However, as show in the case using $\nu^*=0.3$, this behaviour depends on the choice of model parameters and oscillations might appear, indicating a suboptimal adaptive time stepping. Moreover, changing the value of the acceleration constants may also affect the simulation time, the deformed interface, and the pattern structures. In all our experiments, we observe that accepted steps converge in less than 4 iterations in both stages of the RK method.

\paragraph{Example 4: Simulation of the full-3D appendage model.} 
We finally turn to the simulation of the mechanochemical interaction in the 
context of skin appendage using a 3D slab as computational domain, and focusing 
on the reaction model proposed in \eqref{reaction:general}. Model parameters 
are chosen according to Table \ref{table:ex04params}. The domains of interest are the 
blocks $\OmD=(0,50)^3$ and $\OmE = (0,50)^2\times(50,75)$, which we discretise into 
unstructured meshes with 104360 and 60067 tetrahedral elements, respectively. Initial concentrations are set according to steady state solutions obtain from  \eqref{reaction:general}, perturbed with a uniformly distributed random variable with variance $10^{-3}$. The adaptive method described in Algorithm \ref{algo:TRBDF2} is employed with controller parameters as in Table \ref{tab:dtAdaptParam1}, except for $\mathrm{A_{TOL}}=10^{-6}$ and $\mathrm{TOL}_N = \mathrm{R_{TOL}}=10^{-3}$. The process 
is run until $T=300$.

\begin{table}[h]
\begin{center}
\begin{tabular}{l}
\hline
\hline\noalign{\smallskip}
$c_f^\mathrm{D} = 20$, $c_g^\mathrm{D} = 1$, 
$r_1 = 0.25$, $r_0 = 1$, $r_2 = 20$, $r_3 = 5$, $r_4 = 5$, 
$M_{11}^\mathrm{D} = M_{22}^\mathrm{D} = 1$, $M_{33}^\mathrm{D} = 35$, $M_{44}^\mathrm{D} = 10$,\\[1.5ex] 
$c_f^\mathrm{E} = 20$, $c_g^\mathrm{E} = 1$ $r_5 = 5.0$, $r_6 = 20.0$, $r_7 = 5.0$, 
$M_{11}^\mathrm{E} = M_{22}^\mathrm{E} = 1$, $M_{33}^\mathrm{E} = 35$, $M_{44}^\mathrm{E} = 10$,\\[1.5ex] 
$E^* = 1000$, $\nu^*=0.475$, $J^*=K^*=1$, $\alpha^\mathrm{E}=9.5$ \\
\noalign{\smallskip}\hline
\hline
\end{tabular}
\caption{Example 4. Parameters for the mechanochemical 
model using \eqref{reaction:general} and \eqref{2way-coupling}. }\label{table:ex04params}
\end{center}
\end{table}

\begin{figure}[!h]
\begin{center}
\includegraphics[width=0.325\textwidth, clip=true, trim={12.5cm 0cm 3.5cm 0cm}]{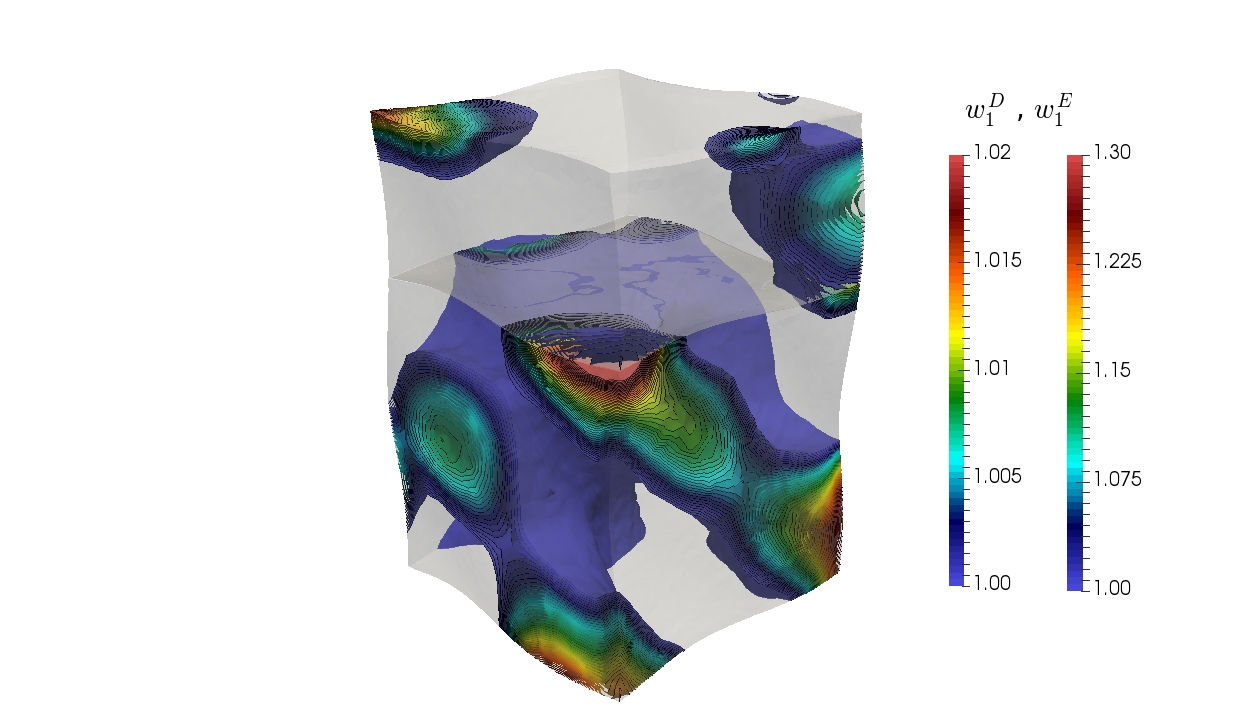}
\includegraphics[width=0.325\textwidth, clip=true, trim={12.5cm 0cm 3.5cm 0cm}]{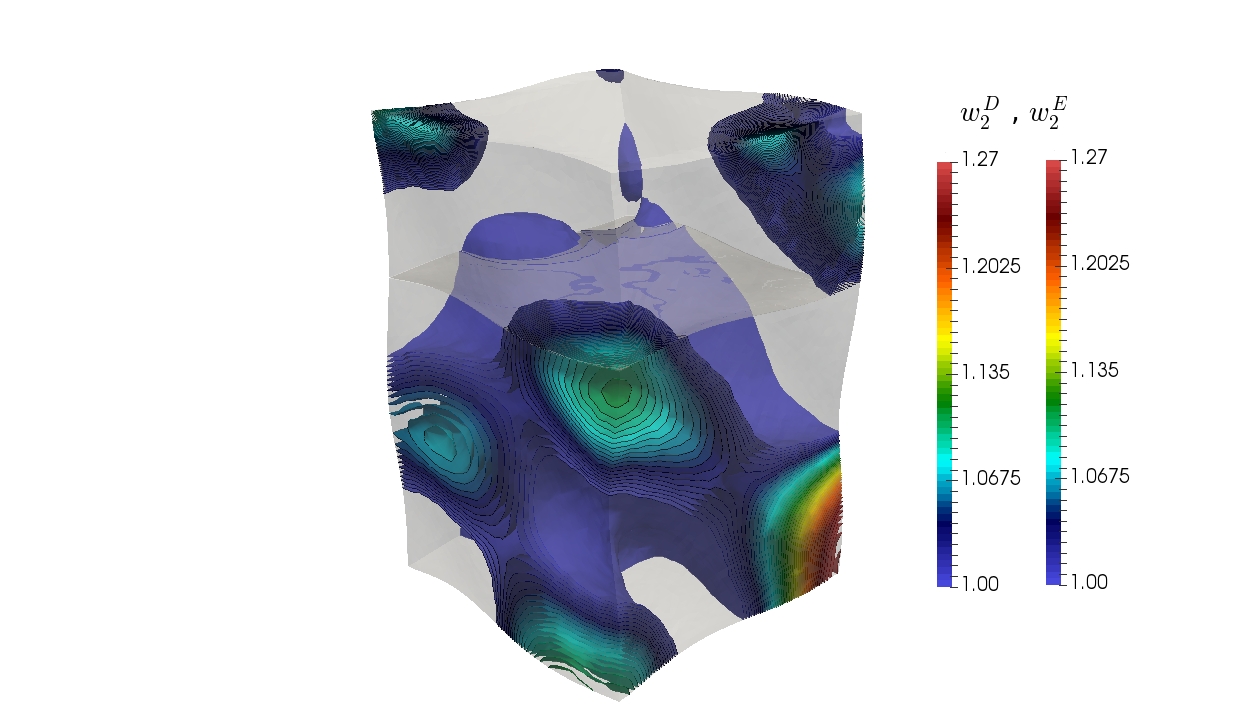}
\includegraphics[width=0.325\textwidth, clip=true, trim={12.5cm 0cm 3.5cm 0cm}]{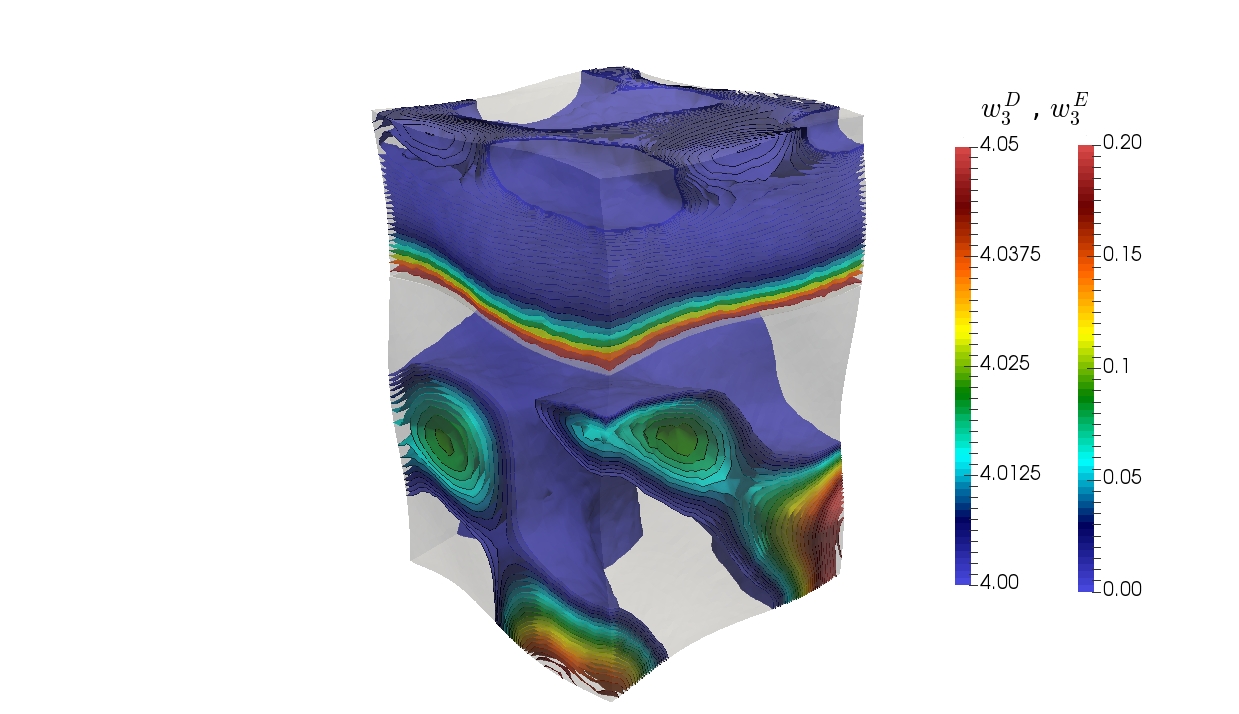}
\includegraphics[width=0.325\textwidth, clip=true, trim={12.5cm 0cm 3.5cm 0cm}]{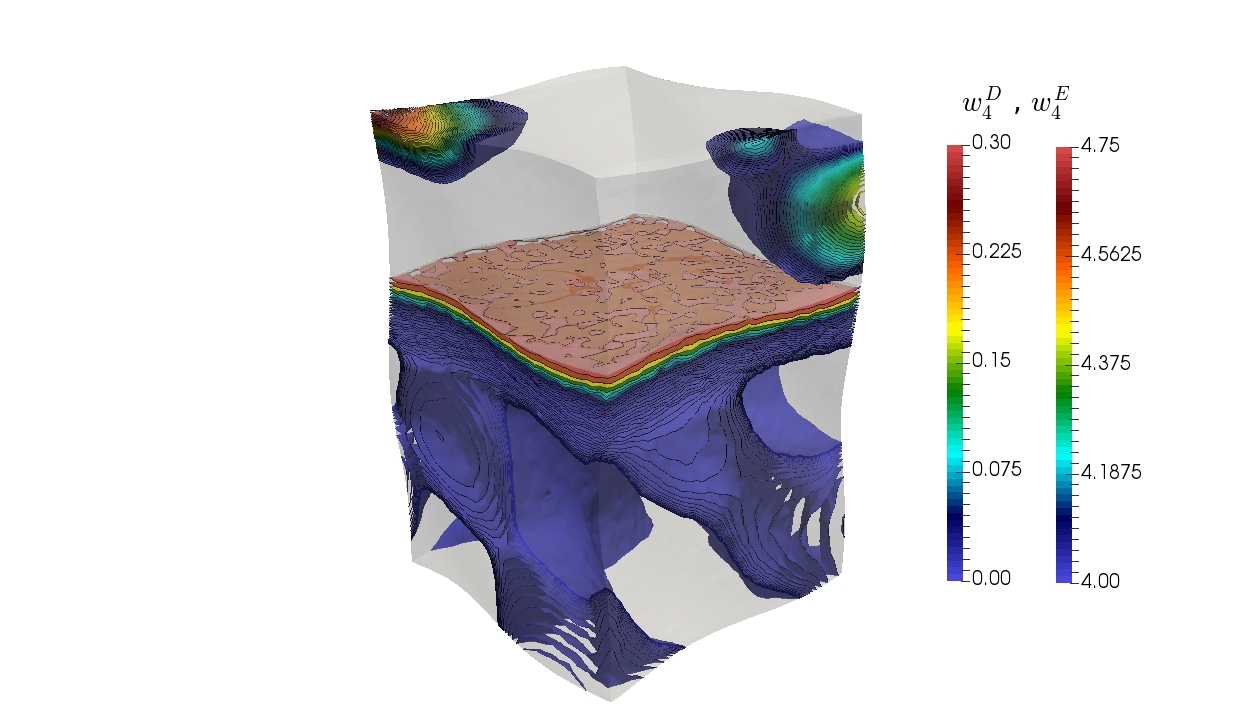}\hspace{1.0cm}
\includegraphics[width=0.325\textwidth, clip=true, trim={12.5cm 0cm 3.5cm 0cm}]{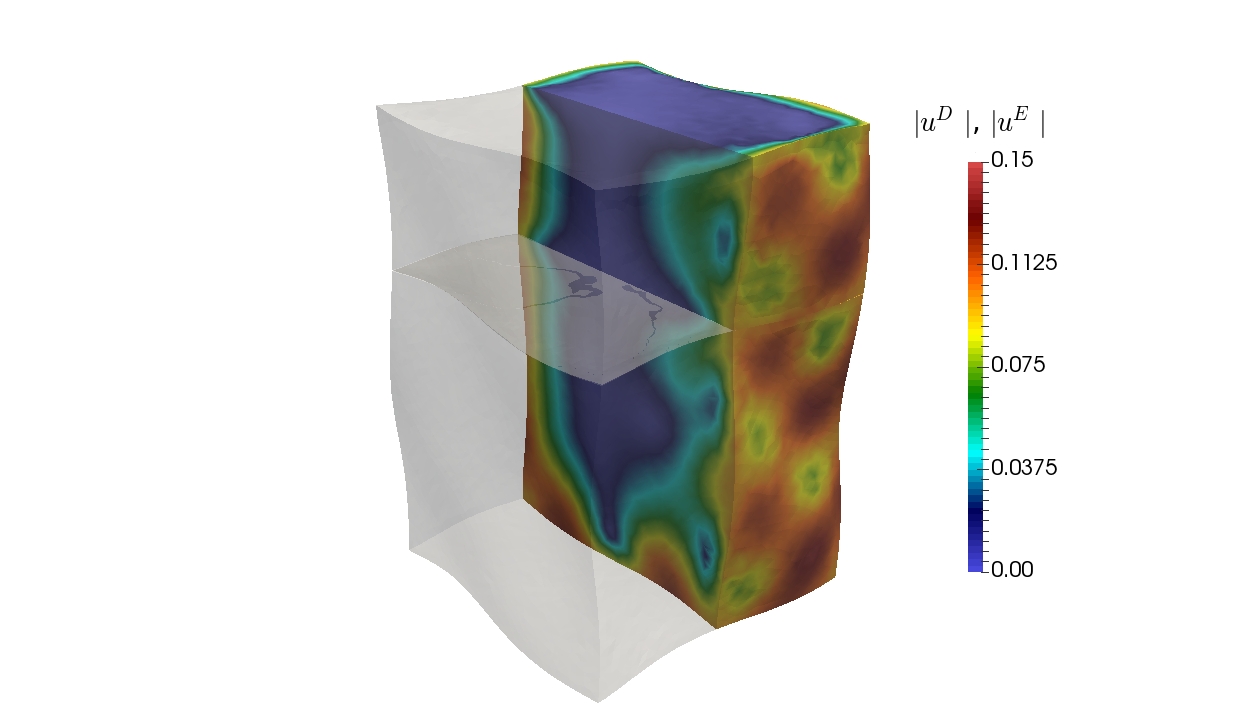}
\caption{Example 4, 3D case. 
Concentration of species $w_1^*$, $w_4^*$ (left and and centre) and displacement magnitude (right), computed a $T=300$ and plotted on the deformed configuration (scaled 10 times for better visualisation). For the displacement magnitude, we cut the domain in half to see the internal displacement.}
\label{fig:ex04-a} 
\end{center}
\end{figure}

The results are presented in Figure \ref{fig:ex04-a}, showing a clear shift in species (cell, matrix and morphogens) concentration patterns across the interface, with the exception of species $w_2^*$ (explained by the similarity in the reaction description \eqref{reaction:general}). The uniform Lam\'e constants on both domains and the Robin conditions imply that deformation patterns are of comparable magnitude throughout $\Omega$, and the zones of high deformations are concentrated near the domain boundary. Different mechanochemical patterns can be produced by changing the coupling and reaction parameters, and a specific study on linear and nonlinear stability will be carried out in a forthcoming contribution. At this stage we can already observe that for $w_1^*$, for example, only the dermis species possesses a production term (with logistic growth) while the dynamics on the epidermis is described by a convection-diffusion process. This explains why the concentration variation in the dermis is much smaller than in the epidermis layer. For $w_3^*$ and $w_4^*$ we observe a similar behaviour, but having a larger variation on the $z$-axis for a specific layer ($\OmD$ for $w_4$, $\OmE$ for $w_3$). This is due to some chemical species being produced only on a specific layer (see \eqref{reaction:general}) and the relevant concentration decreases rapidly, due to metabolism degradation, as one moves away from that layer.

\section{Concluding remarks}\label{sec:concl}

We have introduced a mechanochemical model for the simulation of basic 
processes related to skin appendage patterning mechanisms. The two-way coupling 
between linear elasticity and advection-diffusion-reaction systems is achieved 
through a mechanical term depending on the gradient of the chemical species, 
and a pressure-dependent source term representing production of species.  
A weak formulation in a multidomain setting was introduced, and we presented a partitioned fixed-point and 
Schwarz algorithm to decouple the system into four principal blocks (the two domains, dermis and epidermis, being treated for advection-diffusion-reaction and mechanics). The space discretisation uses MINI-elements for the approximation of the displacement-pressure pair, and piecewise linear and continuous Lagrange elements for the species concentrations. A trapezoidal BDF2 method combined with adaptive Runge-Kutta schemes is used as a time advancing 
strategy, and algorithmic details are provided. The convergence properties of the 
proposed methods were studied  in detail, and a set of numerical tests addresses pattern 
generation and the deformation of the interface. 

Anomalous diffusivity (e.g., fractional, stress- or strain-dependent diffusion \cite{cherubini17}) 
was not considered here. Further work is also needed on the detailed linear stability analysis associated to the coupling mechanisms for the specific models at hand, as well as on the application of the developed methods in the simulation of more complex phenomena. Future work will also address extension of our model to non-linear elasticity and viscoelasticity and theoretical considerations on the convergence properties of partitioned and monolithic schemes. 


\end{document}